\documentclass[a4paper,12pt]{article}

\usepackage{lineno}
\usepackage{amssymb}
\usepackage{amsfonts,amsbsy,graphics,epsfig}
\usepackage{graphicx}
\usepackage{amsmath}
\usepackage{color}
\usepackage{multirow}
\usepackage{arydshln}
\usepackage{lineno,hyperref}
\usepackage{xurl}

\def\1{1\kern-.20em {\rm l}}
\usepackage{float}
\usepackage{caption}
\usepackage{subcaption}
\usepackage{graphicx}
\usepackage{multirow}
\usepackage{adjustbox}
\usepackage{enumitem}
\usepackage[normalem]{ulem}
\newtheorem{lemma}{Lemma}
\newtheorem{theorem}{Theorem}

\newtheorem{remark}{Remark}

\newtheorem{corollary}{Corollary}

\newcommand{\Ex}{\ensuremath{\mathbb{E}}}

\usepackage{booktabs}
\usepackage{url}

\usepackage[dvipsnames]{xcolor}
\usepackage{wasysym}
\usepackage{natbib}
\usepackage{cleveref}

\usepackage{adjustbox}
\modulolinenumbers[1]

\usepackage{authblk}

\title{\textbf{$k$-nearest neighbors prediction and classification  for spatial data}}
\medskip 
\author[1,2]{Mohamed-Salem Ahmed}
\author[3]{Mamadou Ndiaye}
\author[4]{Mohammed Kadi Attouch}
\author[5]{Sophie Dabo-Niang}
\affil[1]{\small {\textit{Alicante, 23 Rue Paul Dubrule, 59810 Lesquin, France,
			msahmed0390@gmail.com}}}
\affil[2]{\small {\textit{Univ. Lille, CHU Lille, ULR 2694 – METRICS : Evaluation des technologies de santé et des pratiques médicales, F-59000 Lille, France}}}
\affil[3]{\small {\textit{University Cheikh Anta Diop of Dakar \& Senegalese Institute of Agricultural Research, Dakar, Senegal, mamadoundiaye397@gmail.com}}}
	
\affil[4]{\small {\textit{Department of Mathematics, University Sidi Bel Abbes, Sidi Bel Abbes,  Algeria, attou\_kadi@yahoo.fr}}}
	
\affil[5]{\small {\textit{Université Lille, CNRS, UMR 8524-Laboratoire Paul Painlevé, INRIA-MODAL, Lille, F-59000, France, sophie.dabo@univ-lille.fr}}}

\date{ }

\begin{document}

 \maketitle

\begin{center}
	\rule{1\linewidth}{1.5pt}
\end{center}
\noindent	
{\textbf{Abstract}}\\
\noindent This paper proposes  a spatial $k$-nearest neighbor method for nonparametric prediction of real-valued spatial data and supervised classification for categorical spatial data. The proposed method is based on a double nearest neighbor rule which combines two kernels to control the distances between observations and locations. It uses a random bandwidth in order to more appropriately fit the distributions of the covariates. The almost complete convergence with rate of the proposed predictor is established and the almost sure convergence of the supervised classification rule was deduced. Finite sample properties are given for two applications of the $k$-nearest neighbor prediction and classification rule to the soil and the fisheries datasets
\medskip 

\begin{center}
	\rule{1\linewidth}{1.5pt}
\end{center}

\section{Introduction}
\label{Sec1}
Spatial data is popularly used in many fields such as environmental sciences, geophysics, soil science, oceanography, econometrics, epidemiology, forestry, image processing and so on where the data of interest are  collected across space. The literature on spatial and spatio-temporal models are relatively abundant \citep[see][]{cress}.

\noindent Complex issues arise in spatial or space-time analysis, many of which are neither clearly defined nor entirely resolved, but form the basis for current researches. The most fundamental of these is spatial prediction, namely the reconstruction of a phenomenon over its domain from a set of observed values.  
Data dependency is one of the practical considerations that influences the available spatial prediction techniques. In fact,
spatial data are often dependent and a spatial regression or prediction model must be able to handle this aspect.\\
During the first half of the twentieth century, spatial prediction was studied in the scope of geostatistics, commonly known as kriging. The latter is a spatial interpolation method which allows a linear prediction of a stationary Gaussian spatial process based on a parametric spatial covariance function.
Since its apparition, kriging has been widely extended and generalized to several directions. Then, several linear spatial regression or prediction methods have been widely studied in the literature.\\
However a pre-selected parametric model might be restrictive.
 In response to that,  a related stream of literature focused on translating the theory of
spatial parametric inference to semi-parametric or nonparametric context. The first results in this direction are those of \cite{tran90}. Consequently,  nowadays, a dynamic concerns the deployment of nonparametric methods to spatial statistics including prediction methods \citep[see][]{biau}.
\noindent However, the literature on spatial non-parametric estimation techniques is not as extensive as  that of the parametric context case. For an overview on results and applications considering spatially dependent data  for regression estimation, prediction and classification, we  highlight  the following works: 
\cite{hallin2004}, \cite{dab07},  \cite{Wang2009}\cite{menez}, \cite{dab13},   \cite{younso2017}, \cite{durocher2019estimating},   \cite{garcia2020use}, \cite{oufdou2021forecasting}, \cite{shi2021nonparametric}.\\
\noindent Among the nonparametric methods, the $k$-Nearest Neighbor ($k$-NN) method is of interest here. The $k$-NN kernel regression estimator (see \cite{biau2015lectures,li2020nonparametric}) is a weighted average of response variables in the neighborhood of the values of covariates. It has a significant advantage over the classical kernel estimate. The specificity of the $k$-NN estimator lies in the fact that it is flexible to various presence of heterogeneity in used covariates which makes it accountable for the local structure of the data.  This consists of the choice of an appropriate number of neighbors using a random bandwidth adapted to the local structure of the data and permitting to enhance the knowledge on local data  dependency.\\ The use of the $k$-NN method is still new in the case of spatial data. \cite{jli} proposed a regression estimator of spatial data based on the $k$-NN method. They proved an asymptotic normality result of their estimator in the case of multivariate data. \cite{fan2021image} investigated and improved the $k$-NN algorithm for image classification. A spatio-temporal nonparametric method based on the $k$-NN approach was proposed by \cite{priambodo2021spatio} for forecasting traffic conditions taking into account the high relationship between roads of a given traffic flow.\\
\noindent The lack of spatial nonparametric prediction techniques with an explicit general spatial proximity structure motivated this work. We are mainly interested in the asymptotic properties of the nonparametric prediction for geostatistical spatial processes using the $k$-NN method. The originality of the suggested predictor lies in the fact that it depends on two kernels, one of which controls the distance between observations using random bandwidth and the other controls the spatial proximity. A similar idea was presented in  \cite{menez}, \cite{dab14},  \cite{cam}, \cite{garcia2020use}  in the context of the classical kernel prediction problem for multivariate or functional spatial data.\\
The present paper aims to give an explicit general spatial proximity structure into a $k$-NN predictor  taking advantages of  existing works in the context of multivariate spatial  data.
We derive a double nearest neighbors selection method of the classical  $k$-nearest neighbors  \citep[see][]{jli}. Furthermore, we study the asymptotic results of such predictors as well as a related classifier and give finite sample properties towards two real data applications. The proposed empirical study reveals two core insights about the $k$-NN proposed method: (i) the importance of taking into account the spatial information in the nonparametric predictor or classifier; (ii) the $k$-NN is less sensitive to choice of kernel functions and allows more flexibility regarding the covariate distribution than the kernel method of \cite{dab14}.
\\
\noindent The rest of this article is organized as follows: Section~\ref{section1} describes the regression model and the proposed predictor; Section~\ref{section2} gives the asymptotic properties of the  proposed predictor;   Section~\ref{section2bis} introduces  a supervised classification rule and gives  its  consistency and  Section~\ref{section3} describes two applications: one related to the prediction of some heavy metals in the region of Swiss Jura  and the other concerns the prediction of presence of some fish species in west Africa. The results are discussed in Section~\ref{conclusion}. The proofs of the asymptotic properties are given in  the Appendix.

\section{Modelling and constructing the predictor}\label{section1}
\noindent
Let $\{Z_{\mathbf{i}}=(X_{\mathbf{i}},Y_{\mathbf{i}})\in\mathbb{R}^{d}\times \mathbb{R}\;,\; \mathbf{i}\in \mathbb{N}^{N}\}$ $(d\geq 1)$ be a spatial process  defined over some probability space $(\Omega,\mathcal{A},\mathbb{P})$, $N\in\mathbb{N}^{*}$. This process is observable in $\mathcal{I}_{\mathbf{n}}=\{\mathbf{i}\in\mathbb{N}^{N}\; :\; 1\leq i_{r}\leq n_{r}\; r=1,\ldots,N \}$,  $\mathbf{n}=(n_{1},\ldots,n_{N})\in \mathbb{N}^{N} $, and $\hat{\mathbf{n}}=n_{1}\times \ldots \times n_{N}$, we write $\mathbf{n}\rightarrow \infty$ if $\min\{n_{r}\}\rightarrow +\infty $, for some constant $C$, $\displaystyle {n_{k}}/{n_{i}}\leq C, \; \forall\; 1\leq k,i \leq N $.  Let $\lVert\cdot\rVert$ denote the Euclidean norm in $\mathbb{R}^{N}$ or in $\mathbb{R}^{d}$ and $\mathbb{I}(\cdot)$ is the indicator function. We assume that the relationship between  $\left\lbrace X_{\mathbf{i}}, \, \mathbf{i}\in  \mathbb{N}^{N}\right\rbrace $ and  $ \left\lbrace Y_{\mathbf{i}}, \, \mathbf{i}\in  \mathbb{N}^{N} \right\rbrace $ is described by the following model:
\begin{equation}
Y_{\mathbf{i}}=r(X_{\mathbf{i}})+\varepsilon_{\mathbf{i}},\;  \mathbf{i}\in \mathbb{N}^{N},
\end{equation}
where
\begin{equation}
r(\cdot)=E \left( Y_{\mathbf{i}}|X_{\mathbf{i}}= \,\cdot \right).
\label{eq-reg}
\end{equation}
The function $r(\cdot)$ is assumed to be independent of $\mathbf{i}$ while the noise $\left\lbrace \epsilon_{\mathbf{i}},\, \mathbf{i}\in \mathbb{N}^{N} \right\rbrace $ is assumed to be centered and independent of $ \left\lbrace X_{\mathbf{i}},\,\mathbf{i}\in\mathbb{N}^{N}\right\rbrace $. 
\\We are interested to predict  the spatial process $\left\lbrace  Y_{\mathbf{i}}, \mathbf{i}\in \mathbb{N}^{N} \right\rbrace $  only at locations where the covariate process $\left\lbrace  X_{\mathbf{i}}, \mathbf{i}\in \mathbb{N}^{N} \right\rbrace $  is observed. Therefore, the prediction framework developed here is designed to only study  a situation in which the covariate process $\left\{ X_{\mathbf{i}}, \mathbf{i}\in \mathbb{N}^{N} \right\}$ describes an observable auxiliary/external information. This situation is well known in the domain of geostatistic by \textit{Kriging with external drift} \citep{hengl2003comparison}.\\ Let consider that at a site  $\mathbf{s_0}\in \mathcal{I}_{\mathbf{n}}$ we  {observe}  $X_{\mathbf{s_0}}$ and  $\left\lbrace (X_{\mathbf{i}},Y_{\mathbf{i}})_{\mathbf{i}\in \mathcal{O}_{\mathbf{n}}}\right\rbrace $ and aim to predict $Y_{\mathbf{s_0}}$ where $ \mathcal{O}_{\mathbf{n}}$ is the observed spatial set $\mathbf{s_0}\notin \mathcal{O}_{\mathbf{n}}$. The sub-set $ \mathcal{O}_{\mathbf{n}}$  is  contained in $\mathcal{I}_{\mathbf{n}}$ and has  a finite cardinal which tends to $\infty$ as $\mathbf{n}\rightarrow \infty$.
However, in order to directly integrate the structure of the spatial dependence in the following proposed spatial predictor, we assume  that the observations $\left\lbrace  (X_{\mathbf{i}},Y_{\mathbf{i}})_{\mathbf{i}\in \mathcal{O}_{\mathbf{n}}} \right\rbrace $ are  locally identically distributed as in   \cite{dab14}. 
The latter  means that a substantial number of observations of $ \left\lbrace (X_{\mathbf{i}},Y_{\mathbf{i}})_{\mathbf{i}\in \mathcal{O}_{\mathbf{n}}}\right\rbrace $ have  distributions close to that of $ (X_{\mathbf{s_0}}, Y_{\mathbf{s_0}})$. \\ 
Finally, we assume that $\left\lbrace Y_{\mathbf{i}}, \,\mathbf{i}\in \mathbb{N}^{N}\right\rbrace $ is integrable and  $(X_{\mathbf{s_0}}, Y_{\mathbf{s_0}})$ has the same distribution as some pair $(X,Y)$ where $(X, Y)$ and any couple of $\left\lbrace (X_{\mathbf{i}},Y_{\mathbf{i}})_{\mathbf{i}\in \mathcal{O}_{\mathbf{n}}}\right\rbrace$ is assumed to have unknown  continuous densities with respect to the Lebesgue measure. Let $f_{X,Y}$ and $f$ denote the densities of $(X,Y )$ and $X$  respectively.

\noindent Now, a $k$-NN predictor of $ Y_{\mathbf{s_0}} $ may be defined by  using a random  bandwidth  depending on the observations and kernel weights, as follows:
\begin{equation}
\widehat{Y}_{\mathbf{s_0}}
=\frac{ \sum_{ \mathbf{i} \in \mathcal{O}_{\mathbf{n}}}   Y_{\mathbf{i}} K_{1}\left (\frac{X_{\mathbf{s_0}}-X_{\mathbf{i}}}{H_{\mathbf{n},X_\mathbf{s_0}}}\right) K_{2}\left (h_{\mathbf{n},\mathbf{s_0}}^{-1}\left \|\frac{\mathbf{s_0}-\mathbf{i}}{\mathbf{n}}\right\|\right)}
{  \sum_{ \mathbf{i} \in \mathcal{O}_{\mathbf{n}}}K_{1}\left (\frac{X_{\mathbf{s_0}}-X_{\mathbf{i}}}{H_{\mathbf{n},X_\mathbf{s_0}}}\right) K_{2}\left (h_{\mathbf{n},\mathbf{s_0}}^{-1}\left \|\frac{\mathbf{s_0}-\mathbf{i}}{\mathbf{n}}\right\|\right)
},
\label{PredictorChap3}
\end{equation}
if the denominator is not null otherwise the predictor is equal to the empirical mean. Here,
$K_{1}$ and $K_{2}$ are two kernels from  $\mathbb{R}^{d} $ and $\mathbb{R}$ to  $\mathbb{R}_{+}$  respectively, $\displaystyle\frac{\mathbf{i}}{\mathbf{n}}=\left(\displaystyle\frac{i_{1}}{n_{1}},\cdots,\displaystyle\frac{i_{N}}{n_{N}}\right)$, and
$$\displaystyle h_{\mathbf{n},\mathbf{s_0}}=\min \left\{ h \in \mathbb{R}^{*}_{+} \,:\, \sum_{\mathbf{i} \in \mathcal{O}_\mathbf{n}}  \mathbb{I}_{\left\{\left\|\displaystyle\frac{\mathbf{i}-\mathbf{s_0}}{\mathbf{n}}\right\|< h\right\}}=k_{\mathbf{n}}^{'}\right\}$$ and $$\displaystyle H_{\mathbf{n},X_\mathbf{s_0}}=\min \left\{ h \in \mathbb{R}^{*}_{+} \,:\, \sum_{\mathbf{i} \in\mathcal{V}_{\mathbf{s_0}}} \mathbb{I}_{\{\|X_{\mathbf{i}}-X_{\mathbf{s_0}}\|< h\}}=k_\mathbf{n}\right\}$$
where $k_{\mathbf{n}}^{'}$, $k_\mathbf{n}$ are positive integer sequences and $\mathcal{V}_{\mathbf{s_0}}=\{\mathbf{i}\in\mathcal{O}_\mathbf{n},\;  \Vert\frac{\mathbf{i}-\mathbf{s_0}}{\mathbf{n}}\Vert< h_{\mathbf{n},\mathbf{s_0}} \}$.

 \noindent The bandwidth $H_{\mathbf{n},X_\mathbf{s_0}}$ is a positive random variable depending on $X_{\mathbf{s_0}}$ and observations $\{X_{\mathbf{i}},\, \mathbf{i}\in\mathcal{O}_{\mathbf{n}}\}$.

\noindent An advantage of using this predictor compared to the fully kernel method proposed by \cite{dab14} lies in its easy implementation. In fact, it is easier to choose the  smoothing parameters $ k_{\mathbf{n}}^{'}$ and $ k_\mathbf{n} $ which take their values in a discrete subset than the bandwidths used in the following kernel counterpart of (\ref{PredictorChap3})  \citep{dab14}
\begin{equation}
\widehat{Y}_{\mathbf{s_0}}^{'}
=\frac{ \sum_{ \mathbf{i} \in \mathcal{O}_{\mathbf{n}}}   Y_{\mathbf{i}}K_{1}\left (\frac{X_{\mathbf{s_0}}-X_{\mathbf{i}}}{h_{\mathbf{n}}}\right)  K_{2}\left (\rho_{\mathbf{n}}^{-1}\left \|\frac{\mathbf{s_0}-\mathbf{i}}{\mathbf{n}}\right\|\right)
}
{  \sum_{ \mathbf{i} \in \mathcal{O}_{\mathbf{n}}} K_{1}\left (\frac{X_{\mathbf{s_0}}-X_{\mathbf{i}}}{h_{\mathbf{n}}}\right) K_{2}\left (\rho_{\mathbf{n}}^{-1}\left \|\frac{\mathbf{s_0}-\mathbf{i}}{\mathbf{n}}\right\|\right)
},
\end{equation}
where here the bandwidths $h_\mathbf{n}$, $\rho_{\mathbf{n}}$ are non random.

\noindent In addition, the fact that $ H_{\mathbf{n}, X_\mathbf{s_0}} $ depends on $ X_{\mathbf{s_0}}$ is the main advantage of the methodology. It allows the predictor $\widehat{Y}_{\mathbf{s_0}}$ to be adapted to a local structure of the observations, particularly if they are heterogeneous \cite[see][]{bur}.

\section{Main results}\label{section2}
\noindent To account for  spatial dependency, we assume that the process $\{Z_{\mathbf{i}}=(X_{\mathbf{i}},Y_{\mathbf{i}})\in\mathbb{R}^{d}\times \mathbb{R}\;,\; \mathbf{i}\in \mathbb{N}^{N}\}$ satisfies a mixing condition defined as follows: there exists a function $\varphi(x)\searrow 0$ as $x\rightarrow \infty$, such that
\begin{align}
\alpha\left(\sigma\left(S\right),\sigma\left(S^{\prime}\right)\right)
& =  \sup\left\{ \left\vert \mathbb{P}(A \cap B)- \mathbb{P}\left(A\right) \mathbb{P}\left(B\right)\right\vert ,\text{ \ }A\in\sigma\left(S\right),\text{ \ }B\in\sigma\left(S^{\prime}\right)\right\} \nonumber \\
& \leq  \psi\left(\text{Card}(S),\text{Card}(S^{\prime})\right)\varphi\left(\text{dist}(S,S^{\prime})\right),
\label{cddep}
\end{align}
\noindent where $S$ and $S'$ are two finite sets of sites, $\text{Card}(S)$ denotes the cardinality of the set $S$.  $\sigma\left(S\right)=\left\{ Z_{\mathbf{i}},\mathbf{i}\in S\right\} $ and $\sigma\left(S^{\prime}\right)=\left\{ Z_{\mathbf{i}},\mathbf{i}\in S^{\prime}\right\} $ are $\sigma$-fields generated by the $Z_{\mathbf{i}}$'s, $\text{dist}(S,S^{\prime})$ is the Euclidean distance between $S$ and $S'$, and $\psi(\cdot)$ is a positive symmetric function nondecreasing in each variable.

\noindent We recall that the process  is said to be strongly mixing if $\psi(\cdot) = 1$ \cite[see][]{Dou1994}. As usual, we will assume that  $\varphi(i)$ verifies :
\begin{equation}
\varphi(t)\leq C t^{-\theta}, \qquad  \; \theta>0\, ,\, t\in\mathbb{R}_{+}^{*}, C>0, \text{ a constant},\, 
\label{cd3}
\end{equation}
i.e.  $\varphi\left(t\right)$ tends to zero at a polynomial rate. Exponential rate may also be considered (see for instance \cite{Dou1994}) for more details.
The asymptotic results given in the following concern only the polynomial case whereas similar results may be obtained easily for the exponential rate.\\
The asymptotic properties of the $k-$NN predictor are achieved  under the following assumptions. Let  $D$ and  $C$ denote for  a compact subset in $\mathbb{R}^{d}$ and a strictly positive generic constant respectively. 
\begin{itemize}
	\item[(H1)] $f$ and $r(\cdot)$ are Lipschitzian functions defined on $D$. In addition, $\inf_{x \in D}f(x)>0$.
	\item[(H2)] $\displaystyle k_\mathbf{n}^{'}\sim \hat{\mathbf{n}}^{\gamma}$ and $k_{\mathbf{n}}\sim \hat{\mathbf{n}}^{\tilde{\gamma}}$, where $\gamma, \, \tilde{\gamma}\in ]0.5,1[$ and   $\tilde{\gamma} < \gamma$.
	\item[(H3)] The kernel $K_{1}$ is  bounded, of compact support and
	\begin{equation}
	\forall u \in \mathbb{R}^{d}, \;  K_{1}(u)\leq K_{1}(tu) \quad \forall t\in ]0,1[.
	\label{Kcd}
	\end{equation}
	\item[(H4)] $K_{2}$ is a bounded nonnegative function, and there exist constants $C_1$, $C_2$ and $\rho$ such that
	\begin{equation}
	C_{1}\mathbb{I}_{\left\{\parallel t\parallel \leq \rho\right\}}\leq K_{2}(\parallel t\parallel)\leq C_{2}\mathbb{I}_{\left\{\parallel t\parallel \leq \rho\right\}}\, , \qquad  \forall \; t\in \mathbb{R}^N, \; 0<C_{1}\leq C_{2}<\infty, \rho>0.
	\label{K1cd}
	\end{equation}
	\item[(H5)]  The density  $f_{X_{\mathbf{i}}X_{\mathbf{j}}}$ of $\left(X_{\mathbf{i}},X_{\mathbf{j}}\right)$ is bounded in $D$ and
	$\left|f_{X_{\mathbf{i}}X_{\mathbf{j}}}(u,v)-f_{X_\mathbf{i}}(u)f_{X_\mathbf{j}}(v)\right|\leq C $ for all $\mathbf{i}\neq\mathbf{j}$ and $(u,v)\in D\times D$ .
	\item[(H6)] $\displaystyle \forall n,m\in \mathbb{N}\quad \psi(n,m)\leq C\min(n,m)$ and \\ $(1-s(1-\tilde{\gamma}))\theta>N\left((2+s(2-\tilde{\gamma}))d+2s(2+\gamma-\tilde{\gamma})\right)$ where, $ 2<s<\frac{1}{1-\tilde{\gamma}}$
	\item[(H7)] $\displaystyle \forall n,m\in \mathbb{N}\quad \psi(n,m)\leq C(n+m+1)^{\tilde{\beta}}$, $\tilde{\beta}\geq 1 $ and \\ $(1-s(1-\tilde{\gamma}))\theta>N\left(2+(2+s(2-\tilde{\gamma}))d+s(4+2\tilde{\beta}+2\gamma-3\tilde{\gamma})\right)$ where $ 2<s<\frac{1}{1-\tilde{\gamma}}$
	\item[(H8)]
	The densities $f_{\mathbf{i}}$ and $f_{X_\mathbf{i},Y_\mathbf{i}}$ of $X_\mathbf{i}$ and $(X_{\mathbf{i}}, Y_{\mathbf{i}}) $ are such that
	$$\sup_{x\in D, \mathbf{i}\in \mathcal{V}_{\mathbf{s_0}}} \vert f_\mathbf{i}(x) - f(x)\vert = o(1),\; \sup_{x\in D, \mathbf{i}\in \mathcal{V}_{\mathbf{s_0}}} \vert g_\mathbf{i}(x) - g(x)\vert = o(1)\qquad \mbox{as}\quad \mathbf{n}\rightarrow \infty,$$
	with $g_{\mathbf{i}}(x)=\int y f_{X_{\mathbf{i}},Y_{\mathbf{i}}}(x,y)dy.$\\
	The conditional density $f_{Y_\mathbf{i},Y_\mathbf{j} \vert X_\mathbf{i},X_\mathbf{j}} $ of $(Y_\mathbf{i},Y_\mathbf{j})$ given $(X_\mathbf{i},X_\mathbf{j})$ and the conditional
	density $f_{Y_\mathbf{i}\vert X_\mathbf{j}}$ of $Y_\mathbf{i}$ given $X_\mathbf{j}$ exist and
	$$f_{Y_\mathbf{i},Y_\mathbf{j} \vert X_\mathbf{i},X_\mathbf{j}}(y,t\vert u,v)< C \qquad f_{Y_\mathbf{i}\vert X_\mathbf{j}}(y\vert u)< C,$$
	for all $y,t,u,v, \mathbf{i},\mathbf{j};\, (u,v)\in D\times D$. 
\end{itemize}
\begin{remark}\label{rem1}\textbf{}
	\begin{enumerate}
		
		\item In assumption (H1),  $f$ is Lipschitzian. This is used particularly for the bias term (see the proof of condition $(L1)$ in Lemmas ~\ref{lemtech1} and ~\ref{lemtech}) and it allows with  assumption (H8) to specify the rate of convergence in Corollary~\ref{Cor1}.
		\item The condition on $k_\mathbf{n}$ in assumption (H2) extends those of the number of neighbors assumed by \cite{mul} in the context of dependent functional time series. The condition on $k_{\mathbf{n}}^{'}$ is the same as those assumed by \cite{dab14} on the number of neighbors of the site $\mathbf{s_0}$.
		\item Condition (\ref{Kcd}) on the kernel $K_{1}$ is required in the proofs of Lemma \ref{lemtech1} and Lemma \ref{lemtech}, for more details on this kernel, see \cite{col80}.
		\item Hypotheses (H4)-(H8) are useful in nonparametric estimation of non-stationary spatial data, see \cite{dab14} for more details. In particular, H4 is imposed for the sake of simplicity  of proofs. It is satisﬁed, for instance, by several kernels with compact support such as triangular, biweight, triweight, Epanechnikov, Parzen kernels.
		\item The theoretical results are obtained under a mixing condition which is not really useful in practice. However, in the case of Gaussian spatial processes, the mixing properties may be linked to parametric Gaussian covariance functions which can be correlated to the kernel function on the locations, see  \cite{robinson2011asymptotic,dab14} for some examples.
	\end{enumerate}
\end{remark}
The following theorem gives an almost complete (a.co) convergence \citep[for details on this kind of convergence, see Definition A.1, page 228, of ][]{ferraty2006nonparametric} of the predictor.
\begin{theorem}\label{th1}
	Under assumptions (H1)-(H5), (H8) and (H6) or (H7), we have
	\begin{equation}
	\left|\widehat{Y}_{\mathbf{s_0}}-Y_{\mathbf{s_0}}\right|\displaystyle\mathop{\longrightarrow}_{{\bf
			n}\rightarrow\infty}0 \quad   a.co.
	\label{res1}
	\end{equation}
\end{theorem}

\noindent If $r(\cdot)$ is Lipschitzian we can obtain the rate of almost complete convergence  
stated in the following corollary.
\begin{corollary}\label{Cor1}
	Under assumptions (H1)-(H5), (H8) and (H6) or (H7), as ${\bf
		n}\rightarrow\infty$,
	\begin{equation}
	\left|\widehat{Y}_{\mathbf{s_0}}-Y_{\mathbf{s_0}}\right|= \mathcal{O}\left(\left(\frac{k_\mathbf{n}}{k^{'}_{\mathbf{n}}}\right)^{1/d}+\left(\frac{\log(\hat{\mathbf{n}})}{k_\mathbf{n}}\right)^{1/2}\right)\quad a.co.
	\label{res2}
	\end{equation}
\end{corollary}
\noindent
The results of Theorem  \ref{th1} and  Corollary \ref{Cor1} can be proved easily from the asymptotic results (stated respectively in Lemmas \ref{th1_reg} and \ref{th2_reg}) of the following regression function estimate
$$ r_{\mathrm{kNN}}(x)=\left \{ \begin{array}{cl}
\displaystyle \frac{g_{\mathbf{n}}(x)}{f_{\mathbf{n}}(x)} & \mathrm{if}\: f_{\mathbf{n}}(x) \neq 0;\\
\displaystyle  \overline{Y},& \text{the empirical mean, otherwise},
\end{array}\right.
$$
with
$$
g_{\mathbf{n}}(x)=\frac{1}{\hat{\mathbf{n}}h_{\mathbf{n},\mathbf{s_0}}^{N}H_{\mathbf{n},x}^{d}}\, \sum_{\mathbf{i}\in \mathcal{I}_{\mathbf{n}}, \mathbf{i}\neq \mathbf{s_0}}K_{1}\left (\frac{x-X_{\mathbf{i}}}{H_{\mathbf{n},x}}\right)K_{2}\left (h_{\mathbf{n},\mathbf{s_0}}^{-1}\left \|\frac{\mathbf{s_0}-\mathbf{i}}{\mathbf{n}}\right\|\right)Y_{\mathbf{i}}
$$
$$
f_{\mathbf{n}}(x)=\frac{1}{\hat{\mathbf{n}}h_{\mathbf{n},\mathbf{s}_{0}}^{N}H_{\mathbf{n},x}^{d}}\, \sum_{\mathbf{i}\in \mathcal{I}_{\mathbf{n}},\mathbf{i}\neq \mathbf{s_0}}K_{1}\left(\frac{x-X_{\mathbf{i}}}{H_{\mathbf{n},x}}\right)K_{2}\left(h_{\mathbf{n},\mathbf{s_0}}^{-1}\left \|\frac{\mathbf{s_0}-\mathbf{i}}{\mathbf{n}}\right\|\right)
$$
\begin{lemma}\label{th1_reg}
	Under assumptions (H1)-(H5), (H8) and (H6) or (H7), we have
	\begin{equation}
	\sup_{x \in D}\left|r_{\mathrm{kNN}}(x)-r(x)\right|\displaystyle\mathop{\longrightarrow}_{{\bf
			n}\rightarrow\infty}0 \quad   a.co.
	\end{equation}
\end{lemma}
\begin{lemma}\label{th2_reg}
	Under assumptions (H1)-(H5), (H8) and (H6) or (H7) and if $r(\cdot)$ is Lipschitzian, as ${\bf
		n}\rightarrow\infty$, we have
	\begin{equation}
	\sup_{x \in D} \left \vert r_{\mathrm{kNN}}(x)-r(x)\right \vert= \mathcal{O}\left(\left(\frac{k_\mathbf{n}}{k^{'}_{\mathbf{n}}}\right)^{1/d}+\left(\frac{\log(\hat{\mathbf{n}})}{k_\mathbf{n}}\right)^{1/2}\right)\quad a.co.
	\end{equation}
\end{lemma}
\noindent The  main difficulty in the proofs of these lemmas comes from the randomness of the window $H_{\mathbf{n},x}$. Then, we do not have in the numerator and denominator of $r_{\mathrm{kNN}}(x)$ sums of identically distributed variables. The idea is to frame sensibly $H_{\mathbf{n},x}$ by two non-random bandwidths. Since the proofs of Theorem \ref{th1} and  Corollary \ref{Cor1} come directly from that of Lemmas \ref{th1_reg} and \ref{th2_reg}, they will be omitted.
\\ In the following, we apply the proposed prediction method to supervised classification.

\section{Construction of a $k$-NN supervised classification rule}\label{section2bis}
The aim here is about predicting the unknown nature of
an object, a discrete quantity for example one or zero.  Let an observation
of an object be a $d$-dimensional vector $X$. The unknown
nature of the object is called a \textit{class} and is denoted by $Y$ which
takes values in a finite set $\{1,\ldots ,M\}$. In classification, one
constructs a function (\textit{classifier}) $g$ taking values in $\{1,\ldots ,M\}$ which
represents one's guess $g(X)$ of $Y$ given $X$. 
\newline
We aim to predict the class $Y$
from $X$ at a given location using a sample of this pair of variables at some observed locations. As in the Section \ref{section1}, we assume that the prediction site is $\mathbf{s_0}\in \mathcal{I}_{\mathbf{n}}$ and $(X_{\mathbf{s_0}},Y_{\mathbf{s_0}})$ has the same distribution as $(X,Y)$ and the observations $(X_{\mathbf{i}},Y_{\mathbf{i}})_{\mathbf{i} \mathcal{O}_{\mathbf{n}}}\in $  are locally identically distributed.\\ The mapping $g$ is then defined on $\mathbb{R}^d$ and takes values in $\{1,\ldots ,M\}$. One makes a mistake  on $Y$ if $g(X)\neq Y$, and the probability of error for a classifier $g$ is given by:
\begin{equation*}
L(g)=P\{g(X)\neq Y\}.
\end{equation*}%
It is well known that the Bayes classifier defined by,
\begin{equation*}
g^{\ast }=\underset{g:\mathcal{X}\rightarrow \{1,\ldots ,M\}}{\arg \min }%
P\{g(X)\neq Y\},
\end{equation*}%
is the best possible classifier with respect to the quadratic loss. The minimum
probability of error is called the \textit{Bayes error} and is denoted by $%
L^{\ast }=L(g^{\ast }).$ Note that $g^{\ast }$ depends on the distribution
of $(X,Y)$ which is unknown.\\ An estimator $g_{\mathbf{n}}$ of $g$ is based on the observations $\left\lbrace (X_{\mathbf{i}},Y_{\mathbf{i}})_{\mathbf{i}\in \mathcal{O}_{\mathbf{n}}}\right\rbrace $; $Y$ is predicted by\\ $g_{\mathbf{n}}\left(
X;\,(X_{\mathbf{i}},Y_{\mathbf{i}})_{\mathbf{i}\in \mathcal{O}_{\mathbf{n}}}\right) $. The performance of $g_{\mathbf{n}}$
is measured by the conditional probability of error%
\begin{equation*}
L_{\mathbf{n}}=L\left( g_{\mathbf{n}}\right) =P\left\{ g_{\mathbf{n}}\left(
X;\,(X_{\mathbf{i}},Y_{\mathbf{i}})_{\mathbf{i}\in \mathcal{O}_{\mathbf{n}}}\right) \neq Y\right\} \geq L^{\ast }.
\end{equation*}%
The sequence $\left\{ g_{\mathbf{n}},\,\ \ \mathbf{n}\in \mathbb{N}^{*N}\right\} $ is the  \textit{%
	discrimination rule}. This has been investigated extensively in the literature particularly for independent or time-series data \citep{paredes2006learning, devroye1994strong, devroye19828, hastie1996discriminant}, see the monograph of \cite{biau2015lectures} for more details. \cite{younso2017} has addressed a discrimination kernel rule for multivariate strictly stationary spatial process $\left\lbrace X_{\mathbf{i}}\in \mathbb{R}^d\right\rbrace_{\mathbf{i}\in \mathbb{N}^N}$ and binary spatial classes $\left\lbrace Y_{\mathbf{i}}\in
(0,1)\right\rbrace_{\mathbf{i}\in \mathbb{N}^{N}}$. To the best of our knowledge this last work is the first one dealing with spatial data. \\
In this section, we extend the previous $k$-NN predictor (\ref{PredictorChap3}) in the setting where $Y$ belongs
to $\{1,\ldots ,M\}$.\\
The Bayes classifier $g^{\ast }$ can be approximated by the rule $\left\lbrace g_{\mathbf{n}},\,\ \ \mathbf{n}\geq 1\right\rbrace$  based on the $k$-NN regression estimate $r_{\mathrm{kNN}}(.)$ and defined as%
\begin{equation}
\sum_{ \mathbf{i} \in \mathcal{O}_{\mathbf{n}}}W_{\mathbf{n}\mathbf{i}}(\mathbf{s_0})\hbox{\rm 1\hskip -3pt I}_{\{Y_{\mathbf{i}}=g_{\mathbf{n}}(\mathbf{s_0})\}}=%
\max_{1\leq j\leq M}\sum_{ \mathbf{i} \in \mathcal{O}_{\mathbf{n}}}W_{\mathbf{n i}}(\mathbf{s_0})\hbox{\rm 1\hskip -3pt I}%
_{\{Y_{\mathbf{i}}=j\}}, \label{e14}
\end{equation}%
where
\begin{equation}
W_{\mathbf{n}\mathbf{i}}(\mathbf{s_0})
=\frac{ K_{1}\left (\frac{X_{\mathbf{s_0}}-X_{\mathbf{i}}}{H_{\mathbf{n},X_\mathbf{s_0}}}\right) K_{2}\left (h_{\mathbf{n},\mathbf{s_0}}^{-1}\left \|\frac{\mathbf{s_0}-\mathbf{i}}{\mathbf{n}}\right\|\right)}
{  \sum_{ \mathbf{i} \in \mathcal{O}_{\mathbf{n}}}K_{1}\left (\frac{X_{\mathbf{s_0}}-X_{\mathbf{i}}}{H_{\mathbf{n},X_\mathbf{s_0}}}\right) K_{2}\left (h_{\mathbf{n},\mathbf{s_0}}^{-1}\left \|\frac{\mathbf{s_0}-\mathbf{i}}{\mathbf{n}}\right\|\right)
}.
\end{equation}
Such classifier $g_{\mathbf{n}}$ (not necessarily uniquely determined) is called an
\textit{approximate Bayes classifier.}\\ Let us say that a rule is good if it is \textit{consistent},
that is if, $L_{\mathbf{n}}\rightarrow L$ in probability or almost surely as $%
\mathbf{n}\rightarrow \infty .$ \citep{devroye1994strong}\\ The almost sure convergence of the proposed rule is established in the
following theorem.

\begin{theorem}
	\label{as-conv-discr} If assumptions (H1)-(H5), (H8) and (H6) or (H7), as ${\bf
		n}\rightarrow\infty$, hold  then,
	$$L_{\mathbf{n}}-L^{\ast}\displaystyle\mathop{\longrightarrow}_{{\bf
			n}\rightarrow\infty}0 \qquad   \mbox{almost
		surely}.$$
\end{theorem}
\noindent
The proof of this theorem consists to show that \cite[by Theorem 2.3 in ][]{gyorfi1996probabilistic}
\begin{equation*}
\int_{D}\left\vert r(x)- r_{\mathrm{kNN}}(x)\right\vert f(x)dx \mathop{\longrightarrow}_{{\bf
		n}\rightarrow\infty}0, \qquad   \mbox{almost
	surely.}
\end{equation*}
This last comes directly from Lemma \ref{th1_reg} and the integratibility of the density function.
\section{Numerical results}\label{section3}
\noindent
After uncovering theoretical properties of the proposed methodology, finite sample performances of the proposed $k$-NN predictor and classifier are provided using two real case studies.
The first case concerns concentration prediction of heavy metal content in the soil of the Swiss Jura region. {The} second case study investigates the prediction of the presence of three fish species in {West} Africa which is a problem of particular economic interest in this region.
The presented results support numerical assessments of the $k-$NN prediction and classification methodologies and a
comparative study with kernel and cokriging approaches. 

\subsection{Environmental case study}\label{app_env}
\begin{figure}[h!]
	\centering
	\includegraphics[width=0.4\textwidth]{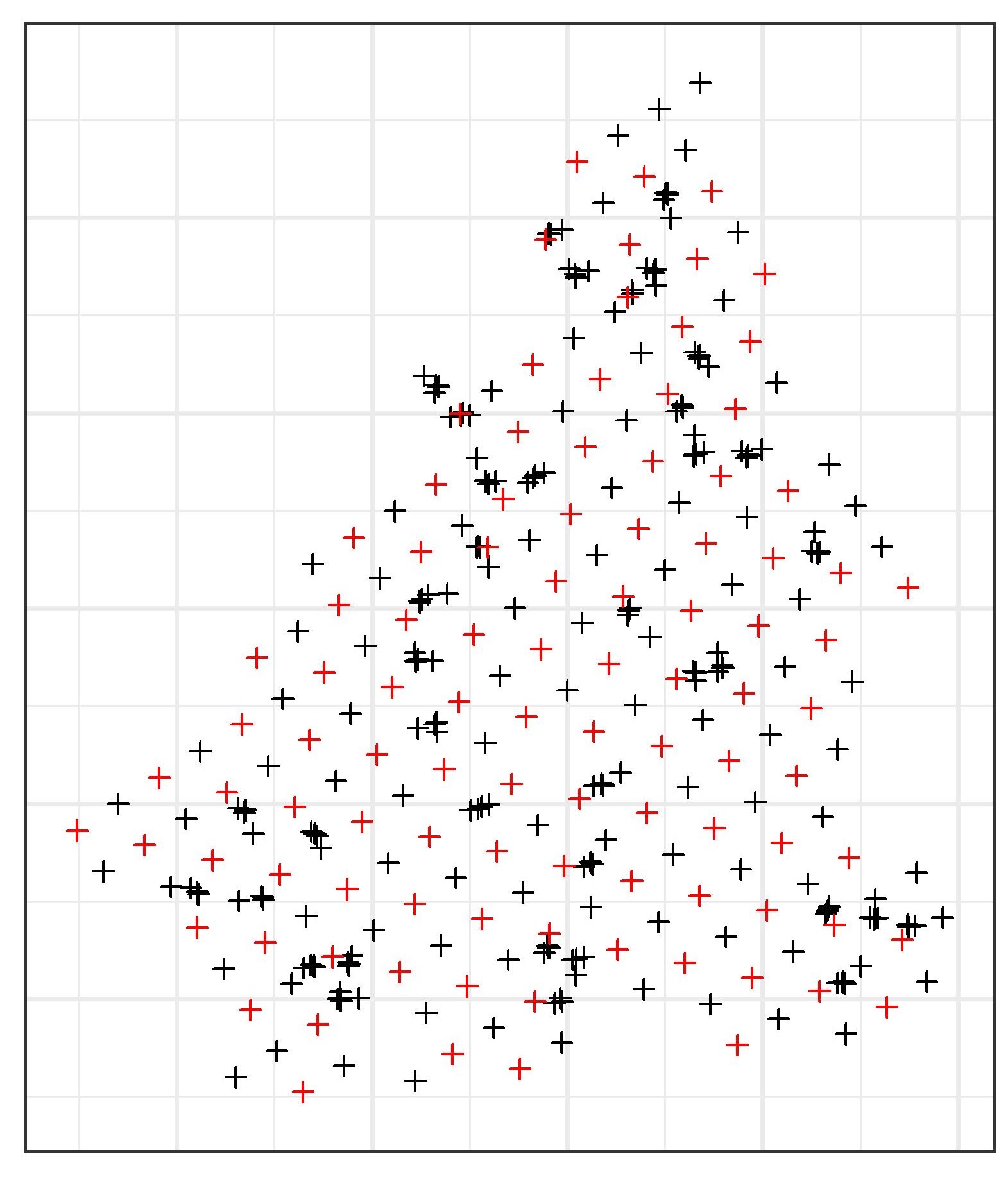}
	\caption{Spatial coordinates in a local grid of $14.5\, \mathrm{Km}^{2}$ in the Swiss Jura study region. The black points denote the training sample, while the red points correspond to the validation sample.  }
	\label{rg}
\end{figure}

\noindent
In this part, we investigate  the performance of the proposed $k$-NN prediction method using the famous Swiss Jura data set ({\url{https://sites.google.com/site/goovaertspierre/pierregoovaertswebsite/download/jura-data}}). This dataset was collected by the Swiss Federal Institute of Technology at Lausanne and studied by several authors \citep[see][e.g]{atteia,goov} in the context of spatial parametric prediction (kriging).
It concerns seven potentially toxic metals (Cadmium, Cobalt, Chromium, Copper, Nickel, Lead and Zinc) measured at $359$ locations represented here by X-Y coordinates in km units in a local grid (Figure \ref{rg}) of $14.5\, \mathrm{km}^{2}$ region in the Swiss Jura.   These locations are divided into two subsets: a training sample which contains $259$ locations and a validation sample which contains $100$ locations, see  Figure \ref{rg}. Refer to \cite{atteia} for more details about the field sampling and laboratory procedures.

\begin{figure}
	\centering
	\includegraphics[width=1\textwidth]{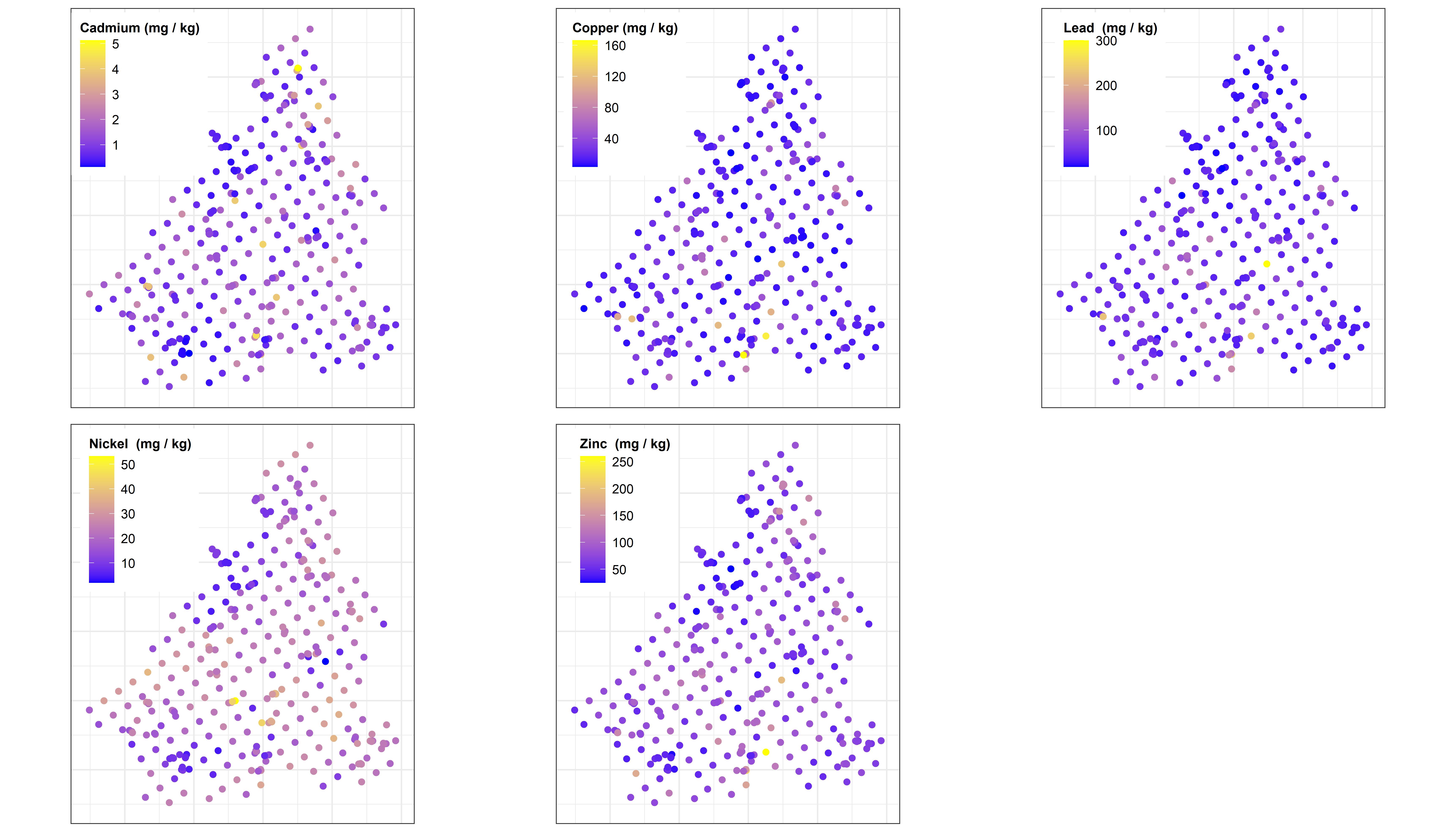}
	\caption{Spatial variation of the concentrations of the five considered toxic metals in the $3$ cases studied.}
	\label{metaux}
	
\end{figure}

\noindent
Many works investigated the performances of many spatial parametric and nonparametric prediction methods through  this dataset. For instance, \cite{goov} considered cokriging  by predicting the concentration in cadmium, copper, and lead at the $100$ validation locations. This study used the secondary variables associated with each case presented in Table~\ref{knn_variables} as covariates. \cite{dab14} investigated the performance of nonparametric (kernel) prediction methods for predicting the concentration of the same three metals and covariates at the $100$ validation sample locations.
We are interested here to compare the performances of our $k$-NN predictor with the kernel predictor of \cite{dab14} and the cokriging one of \cite{goov} using the same training and validation samples as the latter authors. We consider the three cases given in Table \ref{knn_variables} where in each case the response variable $Y_{\mathbf{i}}$'s are given by the observations of the primary variables whereas the covariables $X_{\mathbf{i}}$'s are given by observations of the secondary variables. Figure~\ref{metaux} illustrates the spatial variation  of the concentrations of the five considered metals given in Table~\ref{knn_variables}.  Note that as considered in \cite{goov} and \cite{dab14}, we consider the situation where the secondary variables are assumed to be available at the all $359$ locations. The performances are assessed using the mean absolute error (MAE) of prediction over the $100$ validation sample locations where only the secondary variables are assumed observed. The numbers of neighboring locations $k_\mathbf{n}^{'}$ and observations $k_\mathbf{n}$ used in the $k$-NN prediction method are selected by the leave-one-out cross-validation method. \\ Table~\ref{KNNtb2} gives the MAE of the predictions obtained by the $k-$NN, kernel \citep[those of][]{dab14} and  cokriging \citep[given in][]{goov} methods according to the same combinations of kernels $K_1$ and
$K_2$ which were considered in \cite{dab14}. We use the Euclidean distance on the X-Y coordinates (in km units) defined in the spatial local grid (Figure \ref{rg}). On one hand, it is remarked that the $k$-NN method outperforms its parametric counterpart (cokriging).  This is true for most combinations of kernel functions used in the $k$-NN method. On the other hand, the $k$-NN method gives the best MAE ($ 0.41$) for the prediction of Cadmium while the kernel method outperforms when predicting Copper (MAE $= 6.88$) and Lead (MAE $= 10.06$).\\
The conclusion of this empirical study is  that regardless of the different kernels used, 
the $k$-NN approach outperforms in most situations of  combinations  kernels compare to the kernel method, particularly for cases 2 and 3. It should be note that the kernel method has given the best MAE in the two cases 2 and 3 but it seems to be more dependent to kernels than the $k-$NN method.  
Consistently with previous research \citep[e.g][]{bur}, this may reaffirm the conjecture of robustness of the $k$--NN method and its  flexibility to local
structure of the data compare to the classical kernel method.  A second insight from this empirical study is that the $k-$NN method  may be alternative to the kernel method as well as the  cokriging approach in some situations.

\begin{table}[h!]
	\centering
		\caption{Three considered cases}
	\begin{tabular}{c|l l}
		\hline
		Case & Primary variable & Secondary variables \\
		\hline
		1 &Cadmium & Nickel, Zinc \\
		2 &Copper & Lead, Nickel, Zinc \\
		3& Lead & Copper, Nickel, Zinc\\
		\hline
	\end{tabular}
	\label{knn_variables}
\end{table}

\begin{table}
\centering
\caption{The mean absolute  error of prediction for different parametric and nonparametric methods on the three considered cases. $k$-NN and kernel denote for the mean absolute errors of prediction associated the proposed $k$-NN  predictor and the kernel predictor respectively.  The best mean absolute errors are in bold.}
\begin{tabular}{l l c ll c ll c ll}
\hline
                               &   & \multicolumn{9}{c}{CASE} \\[0.1cm]\cline{4-11}
\multicolumn{2}{c}{Kernels}    &   &  \multicolumn{2}{c}{1} &   &  \multicolumn{2}{c}{2} & & \multicolumn{2}{c}{3}\\[0.1cm] 
\cline{1-2} \cline{4-5} \cline{7-8} \cline{10-11} 
$K_1$       &   $K_2$      & &  $k$-NN & Kernel   & & $k$-NN & Kernel    & & $k$-NN & Kernel \\[0.1cm] \hline
Biweight    & Biweight     & &    0.58 & 0.46 &       & 7.68 & 8.89 &       & 10.69 & 12.62 \\
            & Epanechnikov & &    \textbf{0.41} & 0.47 &       & 7.57  & 9.11 &       & 10.70 & 12.89 \\
            & Gaussian     & &    0.42 & 0.46 &       & 7.38 & 8.88 &       & 10.30 & 12.60 \\
            & Parzen       & &    0.54 & 0.46 &       & 7.30 & 8.86 &       & 10.63 & 12.48 \\
            & Triangular   & &    0.41 & 0.47 &       & 7.41 & 9.14 &       & 10.66 & 12.91 \\
            & Triweight    & &    0.53 & 0.49 &       & 7.46 & 9.26 &       & 10.67 & 13.15 \\[0.15cm]
Epanechnikov& Biweight     & &    0.41 & 0.47 &       & 7.58 & 9.25 &       & 11.09 & 13.15 \\
            & Epanechnikov & &    0.41 & 0.47 &       & 7.74 & 9.16 &       & 10.90 & 12.99 \\
            & Gaussian     & &    0.42 & 0.47 &       & 7.64 & 9.16 &       & 10.62 & 12.98 \\
            & Parzen       & &    0.55 & 0.47 &       & 7.58 & 9.18 &       & 10.98 & 13.04 \\
            & Triangular   & &    0.41 & 0.47 &       & 7.62 & 9.18 &       & 11.09 & 13.06 \\
            & Triweight    & &    0.41 & –    &       & 7.66 & 11.32&       & 11.05 & 15.00 \\[0.15cm]
Gaussian    & Biweight     & &    0.41 & 0.44 &       & 9.67 & 7.02 &       & 14.41 & 11.66 \\
            & Epanechnikov & &    0.41 & 0.44 &       & 9.36 & 7.32 &       & 14.47 & 11.83 \\
            & Gaussian     & &    0.43 & 0.45 &       & 9.14 & 7.91 &       & 14.48 & 12.35 \\
            & Parzen       & &    0.54 & 0.44 &       & 9.08 & 7.91 &       & 14.40 & 12.13 \\
            & Triangular   & &    0.41 & 0.45 &       & 9.08 & 8.28 &       & 14.41 & 12.42 \\
            & Triweight    & &    0.53 & 0.44 &       & 9.91 & \textbf{6.88} &       & 14.38 & \textbf{10.06} \\[0.15cm]
Triangular  & Biweight     & &    0.41 & 0.46 &       & 7.70 & 8.90 &       & 10.94 & 12.62 \\
            & Epanechnikov & &    0.41 & 0.46 &       & 7.76 & 8.91 &       & 10.92 & 12.86 \\
            & Gaussian     & &    0.42 & 0.46 &       & 7.63 & 8.86 &       & 10.46 & 12.50 \\
            & Parzen       & &    0.54 & 0.46 &       & 7.51 & 8.90 &       & 10.76 & 12.86 \\
            & Triangular   & &    0.41 & 0.47 &       & 7.61 & 9.14 &       & 10.92 & 12.90 \\
            & Triweight    & &    0.53 & 0.49 &       & 7.64 & 9.26 &       & 10.89 & 13.15 \\[0.15cm]
Triweight   & Biweight     & &    0.42 & 0.50 &       & 7.40 & 9.27 &       & 10.45 & 13.31 \\
            & Epanechnikov & &    0.42 & 0.50 &       & 7.45 & 10.42&       & 10.45 & 14.40 \\
            & Gaussian     & &    0.43 & 0.50 &       & 7.26 & 10.44 &      & 10.12 & 14.54 \\
            & Parzen       & &    0.55 & 0.50 &       & 7.16 & 10.38 &      & 10.25 & 14.39 \\
            & Triangular   & &    0.42 & 0.50 &       & 7.30 & 9.27 &       & 10.39 & 13.20 \\
            & Triweight    & &    0.53 & –    &       & 7.27 & 11.41 &      & 10.46 &15.11 \\[0.15cm]
\multicolumn{11}{l}{\textit{Parametric methods}:}\\
\multicolumn{2}{l}{Ordinary Cokriging}        & & \multicolumn{2}{c}{0.51} && \multicolumn{2}{c}{7.90} && \multicolumn{2}{c}{10.80} \\
\multicolumn{2}{l}{Revisited Cokriging (cov)} && \multicolumn{2}{c}{0.52} && \multicolumn{2}{c}{7.80} && \multicolumn{2}{c}{10.70}  \\
\multicolumn{2}{l}{Revisited Cokriging (corr)} && \multicolumn{2}{c}{0.52} && \multicolumn{2}{c}{7.40} && \multicolumn{2}{c}{10.60}  \\ \hline 
\end{tabular}
\label{KNNtb2}
\end{table}
\subsection{Fisheries case study}
\begin{figure}[!h]
	\centering
		\includegraphics[width=0.235\textwidth]{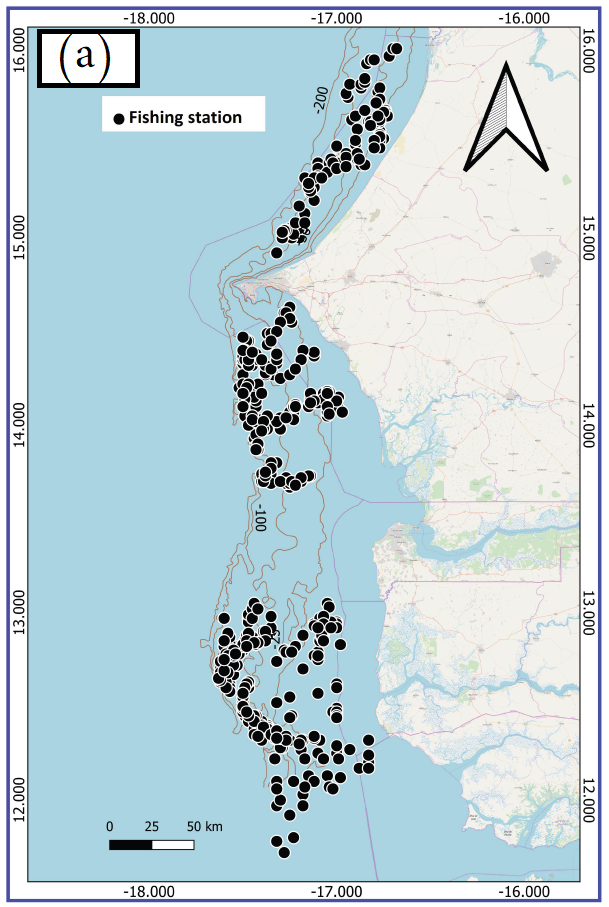}
		\includegraphics[width=0.25\textwidth]{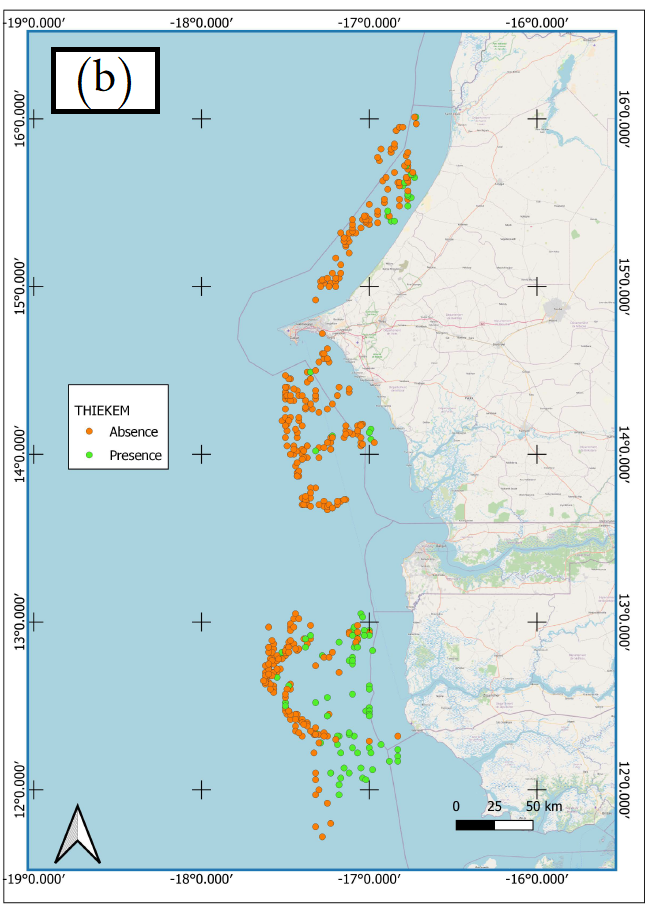}
		\includegraphics[width=0.25\textwidth]{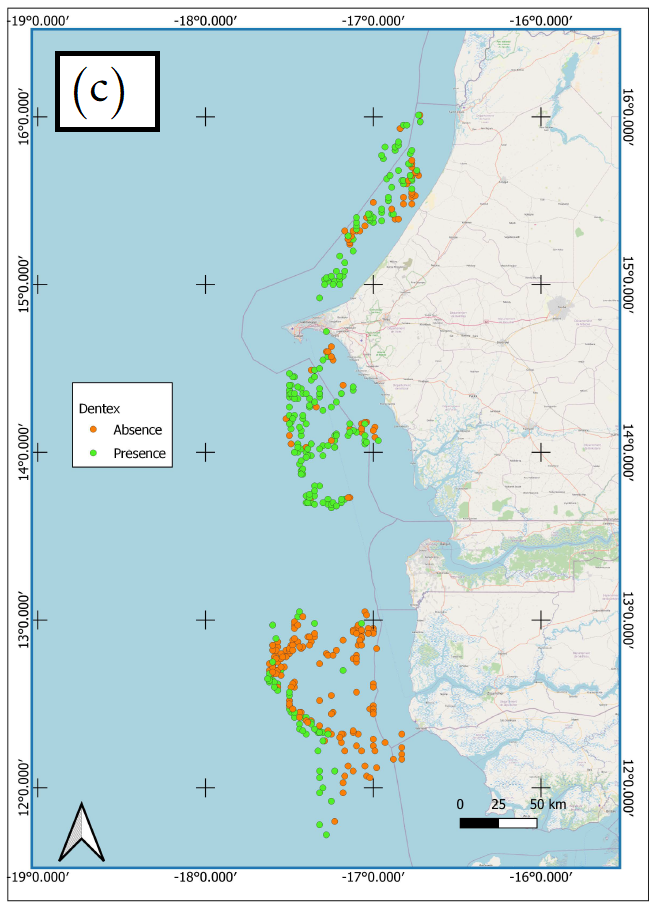}
		\includegraphics[width=0.25\textwidth]{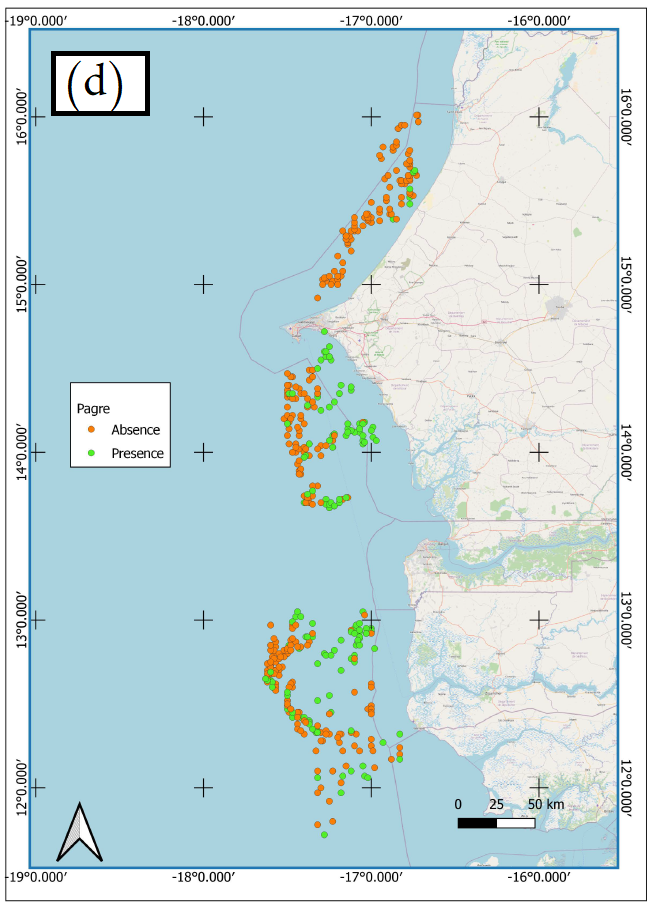}
	\caption{Panel (a) shows the positions of the fishing stations. Panels (b), (c) and (d) show the spatial locations of the three considered species respectively: \textit{Thiekem}, \textit{Dentex}, \textit{Pagrus}. Green points indicate the presence of coastal demersal fish while orange points indicate absence.}
	\label{plotPoint}
\end{figure}

\begin{figure}[!h]
	\centering
		\includegraphics[width=0.25\textwidth]{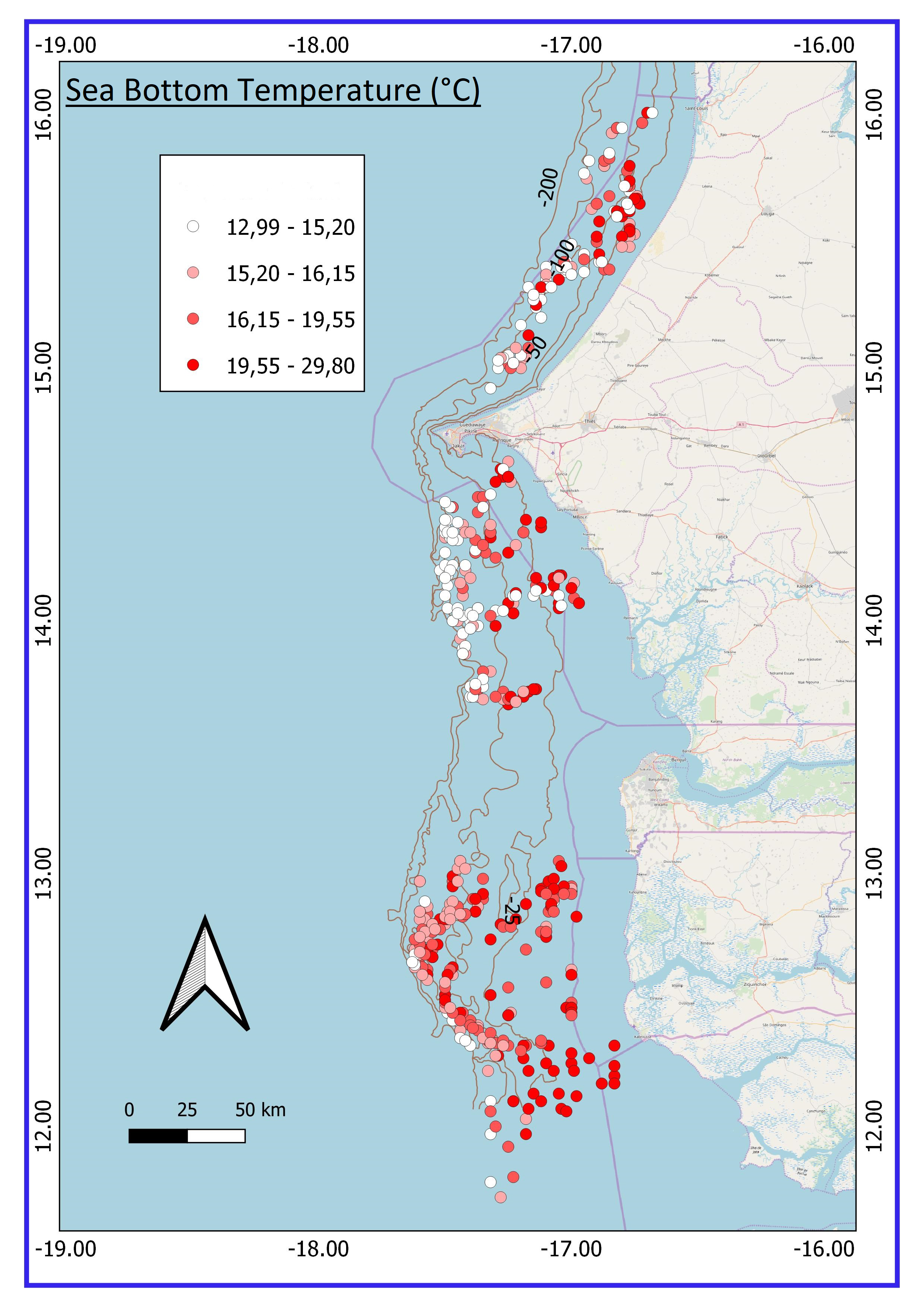}
		\includegraphics[width=0.25\textwidth]{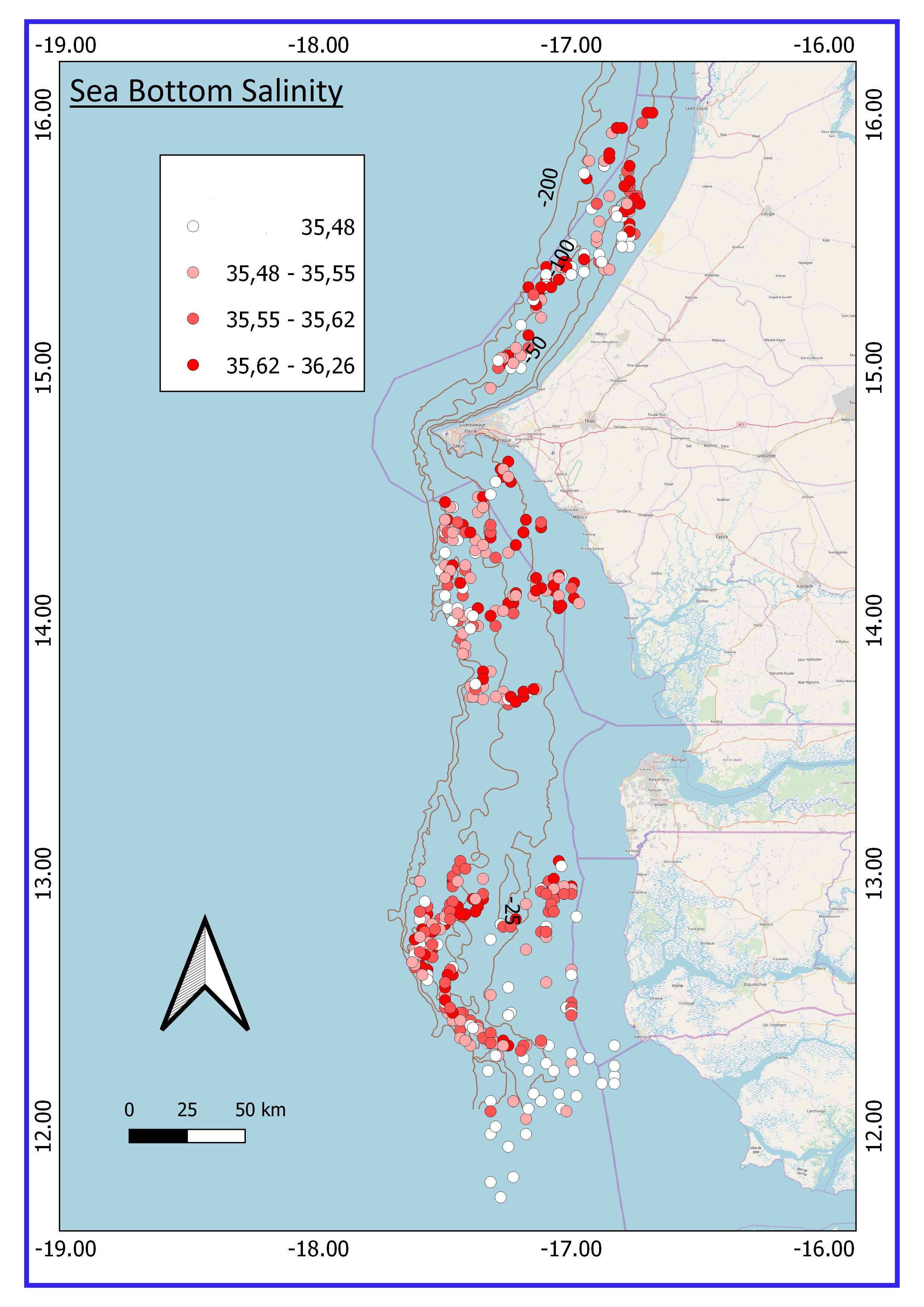}
		\includegraphics[width=0.25\textwidth]{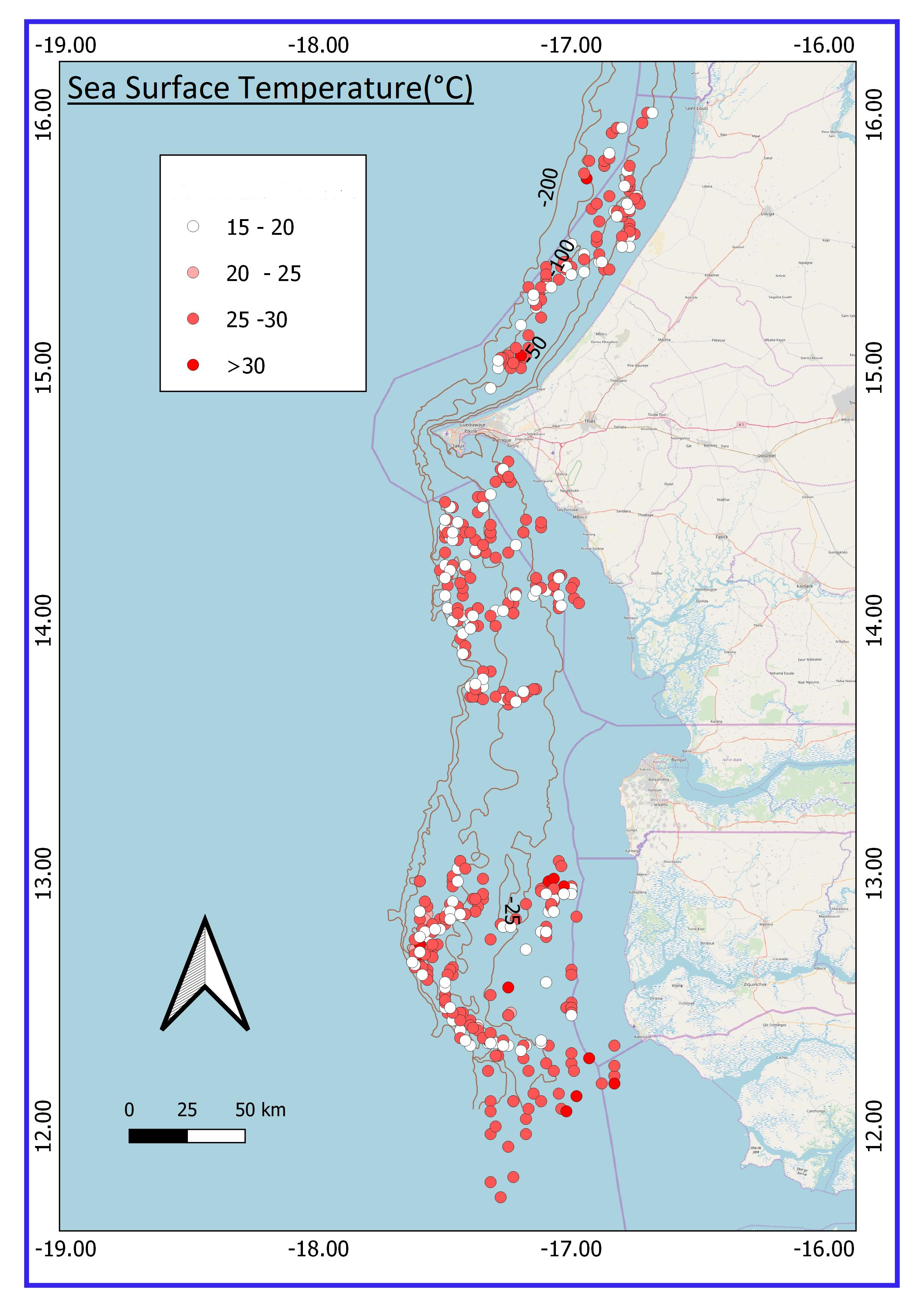}
		\includegraphics[width=0.25\textwidth]{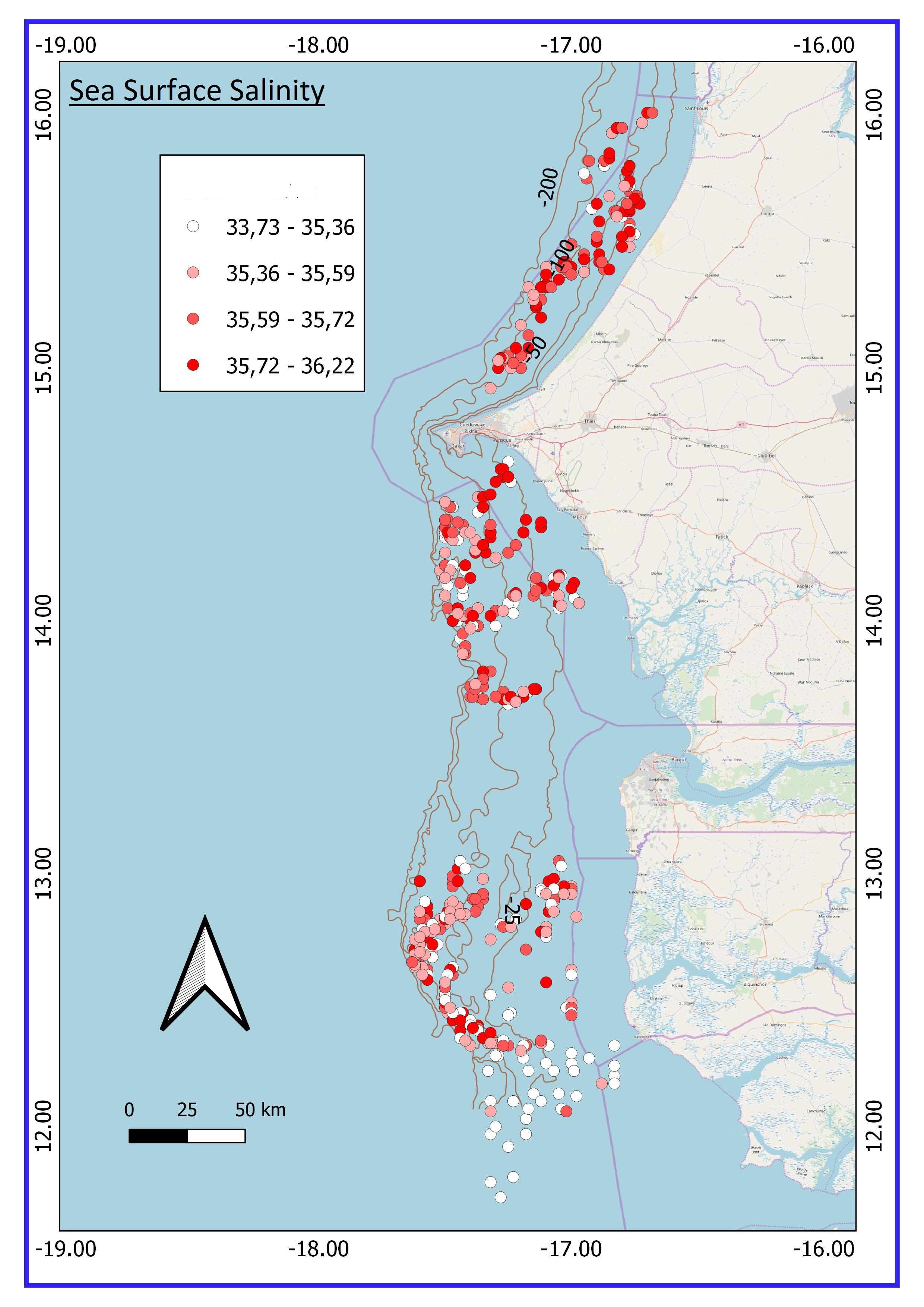}
	\caption{Spatial  variation of temperature and salinity in the bottom and surface.}
	\label{FigImage}
\end{figure}

\begin{figure}[!h]
\centering
\includegraphics[width=1\textwidth]{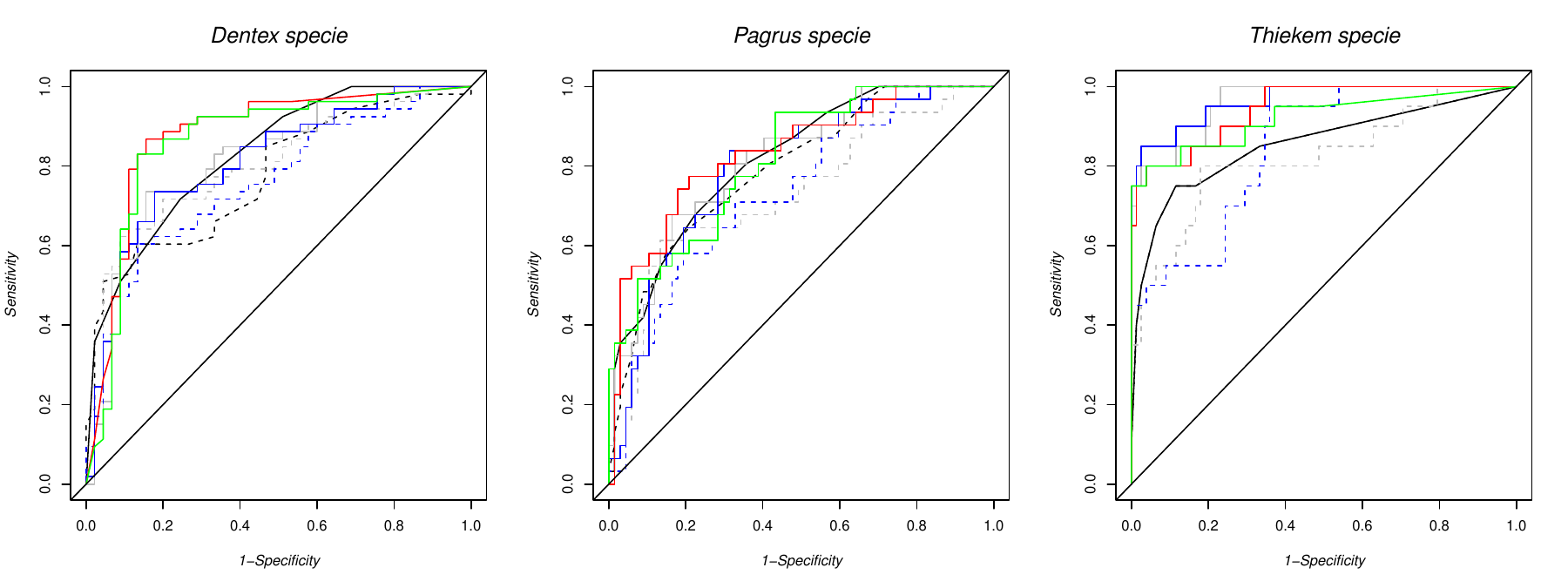}
\includegraphics[width=1\linewidth]{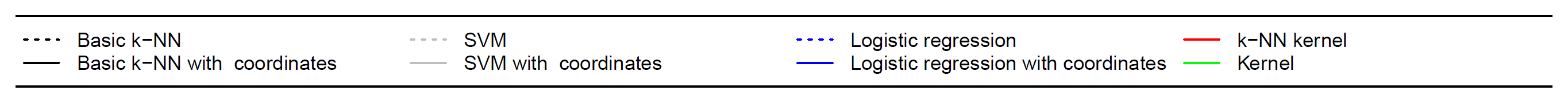}\\
\caption{Comparison of receiver-operating characteristic curves for the various classifiers over the test sets associated with the three considered species respectively.}
\label{fig_fish_appli}
\end{figure}

We consider data from the coastal demersal sea surveys of Senegal performed by the scientific team of the Oceanographic Research Center of Dakar-Thiaroye and  the oceanographic research center of the Senegalese Institute of Agricultural Research, during the cold and hot seasons in the North, Center and South areas of the Senegalese coasts.
Fishing stations were visited from sunrise to  sunset (diurnal strokes) at the rate of half an hour per station. They were selected using stratified sampling, following double stratification by area (North, Center and South) and bathymetry ($0-50~\mathrm{m}$, $50-100~\mathrm{m}$, $100-150~\mathrm{m}$ and $150-200~\mathrm{m}$). The database includes $495$ stations (see panel (a) of Figure \ref{plotPoint}), described among others by  the campaign, temporal features (season, starting and ending trawl times, duration time), spatial coordinates (starting and ending latitude and longitude, area, starting and ending depths, average depth and bathymetric strata), biological parameters (species, family, zoological group and specific status) and environmental parameters (sea bottom temperature  (SBT), sea surface temperature (SST), sea bottom salinity (SBS) and sea surface salinity  (SSS)).\\
It should be noted that  the Senegalese and Mauritanian upwellings affect the spatial and seasonal distributions of coastal demersal fish. Thus, it is important to study the locations of the fish species in this region. In this section, we focus on the three following species which have a particular economic interest in the west African region:
\begin{itemize}
\item[\textbullet] {\it Galeoides decadactylus} (\textit{Thiekem}) of the Polynemidae family that belongs to the coastal community of Sciaenidae. It is  located at a depth between $10~\mathrm{m}$ and $20~\mathrm{m}$. We can say that it is present up to $60~\mathrm{m}$. Panel (b) of Figure \ref{plotPoint} shows the $101$ stations where it is identified. This species is particularly abundant in the south of the Senegal coast. 
\item[\textbullet] {\it Dentex angolensis} (\textit{Dentex}) of the Sparidae family located in tropical and temperate regions. {\it Dentex angolensis} is the most deep species of  Sparidae family. It is present at depths up to $200~\mathrm{m}$. Panel (c) of Figure \ref{plotPoint} shows the $266$ stations where it is identified. This species is particularly abundant in the center and the northern parts of the Senegal coast.
\item[\textbullet] {\it Pagrus caeruleostictus} (\textit{Pagrus}) of the inter-tropical species belonging to the Sparidae family {is} abundant in the south of Dakar (center zone) between $15$ and $35~\mathrm{m}$. It prefers cold waters $(<15 ^{\circ} C)$ between $10$ to $90 ~\mathrm{m}$. Panel (d) of Figure \ref{plotPoint} illustrates the $158$ stations where this species is identified. The \textit{Pagrus} species is mainly present in the center of  the Senegal coast.
\end{itemize}
Figure~\ref{plotPoint} illustrates the spatial distributions of the previous three species. For example, one can observe that {\it Thiekem} is a coastal species. It prefers a higher temperature and lower surface salinity, see Figure~\ref{FigImage} for the spatial distributions of the environmental predictors. One can observe the spatial heterogeneity of the environmental predictors which partially may determine the vertical and horizontal migration of species.
\noindent
We aim to  predict the presence of the three fish species at a given station (location, where only covariates are assumed to be observed) by using the proposed classification rule. We set the response variable $Y$ as $1$ for presence of a species and $0$ otherwise, at each one of $495$ stations. The four environmental variables are considered as the covariates. Our $k$-NN classifier is compared with the kernel classifier derived from the regression estimate proposed by \cite{dab14} and  the following three standard classification methods:
\begin{itemize}
\item[\textbullet] The basic $k$-NN classifier given by the \texttt{caret} package of the \texttt{R} software, where the number of neighbors is chosen by cross validation (CV). We consider two basic $k$--NN classifiers: one which only uses four environmental variables and the second which uses the geographical coordinates (longitude and latitude) in addition to the environmental variables. Note that this second classifier accounts for some spatial proximity. 
\item[\textbullet] SVM (Support vector machines) with radial basis function defined in the \texttt{cadet} package. Analogous to the basic $k$-NN classifier, we consider two SVM, one with only four of the  environmental variables as covariates and the second additionally uses the geographical coordinates.
\item[\textbullet] Logistic regression models with the best model selected by AIC (Akaike information criteria) using a forward-backward variable selection procedure. The latter is applied as before on a model containing only environmental covariates variables and on another containing both environmental variables and geographic coordinates.  
\end{itemize}
In order to compare the different classifiers for each fish species, {the data set is randomly stratified with respect to  the distribution of the outcome variable and the spatial area of the Senegalese coasts (North, Center and South)} into two samples: training and test (validation) samples with respective sizes of $80\%$ and $20\%$ of the total sample size. The training sample was used to construct each classifier using the criteria mentioned above  to select the optimal tuning parameters associated with each classifier. The performance of the different classifiers was compared based on six criteria: area under the receiver operating characteristic curve (AUCOR), accuracy, sensitivity, specificity, negative positive value rate (NPV), and positive predicted value rate (PPV). It should be noted the ROC curve is used to determine the best cutoff point associated with each classifier.  The kernels used in our $k$--NN classifier and the kernel classifier are selected during the parameter setting step from the set of kernels used in Section~\ref{app_env}. Note that for these  two classifiers, the great circle distance was used to calculate the proximities between the spatial locations defined by latitude and longitude coordinates.  \\    
\noindent The results over the test samples of the three fish species are presented in Table~\ref{tab_fish_appli} and Figure~\ref{fig_fish_appli}. Hence, we can remark that incorporating spatial information is of importance. This is evident as the classifiers using spatial information outperform  those which ignore this information.  \\
For the  \textit{Dentex} species, the proposed $k$-NN classifier gives the best AUROC ($88\%$), accuracy ($85\%$), sensitivity ($87\%$) and NPV rate ($84\%$). The basic $k$-NN with spatial coordinates has the best specificity score ($91\%$) while the SVM classifier has the best PPV rate ($89\%$). \\
For the \textit{Pagrus} species, the proposed $k$-NN classifier again gives  the best AUROC ($83\%$), accuracy ($78\%$), sensitivity ($77\%$) and NPV rate ($88\%$). The basic $k$-NN with spatial coordinates has the best specificity ($87\%$) and PPV rate (65\%). Our $k$-NN classifier performs slightly better than the kernel classifier.
\\
For the \textit{Thiekem} species, the logistic regression with spatial coordinates and SVM with spatial coordinates outperform  the others classifiers in terms of AUROC, see Figure~\ref{fig_fish_appli}.  The logistic regression with spatial coordinates gives the best sensitivity ($90\%$) and NPV ($97\%$) while the kernel classifier gives the best accuracy ($92\%$), specificity ($95\%$), and PPV rate ($80\%$).

\begin{table}[!h]
\centering
\caption{Operating characteristics of  the various classifiers over the test samples associated with the three considered fish species. The best metric values are in bold.}
\begin{adjustbox}{max width=\textwidth}
\begin{tabular}{l c c c c c c l l}
\hline 
\multirow{2}{*}{Method}  & \multicolumn{6}{l}{Operating characteristic:}&  \multicolumn{2}{l}{Selected kernels:} \\
& AUROC   & Accuracy & Sensitivity & Specificity & PPV   & NPV   &  $K_1$ &  $K_2$ \\ \hline
\multicolumn{9}{c}{Dentex specie}\\\hline 
Basic k-NN & 0.77 &	0.66 &	0.60	& 0.73 &	0.73 &	0.61 &       &  \\
Basic k-NN with  coordinates & 0.83 &	0.69 &	0.51&	\textbf{0.91}&	0.87&	0.61&       &  \\
SVM   & 0.80   & 0.76  & 0.62  & 0.91  & \textbf{0.89}  & 0.67  &       &  \\
SVM with  coordinates & 0.82  & 0.78  & 0.72  & 0.84  & 0.84  & 0.72  &       &  \\
Logistic regression & 0.76  & 0.70   & 0.60   & 0.82  & 0.80   & 0.64  &       &  \\
Logistic regression with  coordinates & 0.81  & 0.72  & 0.60   & 0.87  & 0.84  & 0.65  &       &  \\
$k-$NN kernel & \textbf{0.88}  & \textbf{0.85}  & \textbf{0.87}  & 0.82  & 0.85  & \textbf{0.84}  & Epanechnikov & Triweight \\
Kernel & 0.86  & 0.84  & 0.81  & 0.87  &  0.88  & 0.80  & Epanechnikov & Triweight \\\hline
\multicolumn{9}{c}{Pagrus specie}\\\hline 
Basic $k-$NN & 0.79	&0.74&	0.65&0.79&0.59&0.83&       &  \\
Basic $k-$NN with coordinates & 0.81&0.77&0.55&\textbf{0.87}	&\textbf{0.65}&	0.81 &       &  \\
SVM   & 0.73  & 0.62  & 0.68  & 0.60     & 0.44  & 0.80   &       &  \\
SVM with  coordinates & 0.81  & 0.70   & 0.74  & 0.69    & 0.52  & 0.85  &       &  \\
Logistic regression & 0.73  & 0.68  & 0.71  & 0.67    & 0.50   & 0.83  &       &  \\
Logistic regression with  coordinates & 0.79  & 0.71  & 0.68  & 0.73    & 0.54  & 0.83  &       &  \\
$k-$NN kernel & \textbf{0.83}  & \textbf{0.78}  & \textbf{0.77}  & 0.78    & 0.62  &\textbf{0.88}  & Gaussian & Triangular \\
Kernel & 0.81   & 0.71  & 0.61  & 0.76    & 0.54   & 0.81  & Biweight & Triangular \\\hline
\multicolumn{9}{c}{Thiekem specie}\\\hline 
Basic k-NN & 0.83 &	0.84&	0.60	&0.9&	0.60	&0.90 &       &  \\
Basic k-NN with  coordinates & 0.86	& 0.88&	0.65& 0.94&0.72 &0.91&       &  \\
SVM   & 0.82  & 0.78  & 0.80   & 0.77  & 0.47  & 0.94  &       &  \\
SVM with  coordinates & \textbf{0.96}  & 0.90   & 0.80   & 0.92  & 0.73 & 0.95  &       &  \\
Logistic regression & 0.85  & 0.72  & 0.70   & 0.73  & 0.40   & 0.90   &       &  \\
Logistic regression with  coordinates & \textbf{0.96}  & 0.84  & \textbf{0.90}   & 0.82  & 0.56  & \textbf{0.97}  &       &  \\
$k$-NN kernel  & 0.94  & 0.91  & 0.80  & 0.94  & 0.76  & 0.95  & Triangular & Triweight \\
Kernel  & 0.92  & \textbf{0.92}  & 0.80   & \textbf{0.95}  & \textbf{0.80}  & 0.95  & Triangular & Parzen \\\hline
\end{tabular}
\end{adjustbox}
\label{tab_fish_appli}
\end{table}

\section{Conclusion}\label{conclusion}
In this work, we proposed a nonparametric spatial $k$-NN prediction for real valued spatial data and derived a supervised classification rule for categorical spatial data. The proposed $k$-NN method combines two kernels to control the distances between  observations and locations as in \cite{dab14}. It uses a random bandwidth defined by the $k$-th lower distance between the covariate prediction point and the covariates of the training sample. The use of a random bandwidth allows more flexibility regarding the covariate distribution. We  have  established infinite sample properties towards the almost complete convergence with rate of the proposed predictor. Then, an almost sure convergence result of the supervised classification rule was deduced. \\
The proposed method has been applied to environmental data prediction and fisheries data classification. The method was used to predict the level of some heavy metals in the Swiss Jura. This application showed that the proposed $k$-NN prediction method outperforms the kernel method of \cite{dab14} and the standard cokriging prediction method. Secondly, the supervised classification rule was applied to predict the presence of three fish species in west Africa. This application is of economic importance in this part of the world. This shows  that  the proposed $k$-NN classifier may be an alternative to the kernel classifier and other well known classifiers (SVM, logistic, basic $k$-NN). Hence, we  argue that the proposed nonparametric prediction method may be a good alternative in some situations compared to the spatial kernel method of \cite{dab14} or usual parametric methods.
\section{Appendix} \label{appendix}
We start by the following technical lemmas that are helpful to handle the difficulties induced by the random bandwidth $H_{\mathbf{n},x}$ in $r_{\mathrm{kNN}}(x)$. They are adaptation of the results given in  \cite{col80} (for independent multivariate data) and their generalized version by \cite{bur}, \cite{Kudra2013} (for independent functional data).
\subsection*{\textbf{Technical Lemmas}}
\noindent
For any random positive variable $T$, $\mathbf{n}\in \mathbb{N}^{*N}$, and $x \in D$, we define

\begin{equation*}c_{\mathbf{n}}(T)=\frac{\sum_{\mathbf{i}\in \mathcal{I}_{\mathbf{n}},\mathbf{s_0}\neq \mathbf{i}}Y_{\mathbf{i}}K_{1}\left(\frac{x-X_{\mathbf{i}}}{T}\right) K_{2}\left(h_{\mathbf{n},\mathbf{s_0}}^{-1}\left \|\frac{\mathbf{s_0}-\mathbf{i}}{\mathbf{n}}\right\|\right)}{\sum_{\mathbf{i}\in \mathcal{I}_{{\mathbf{n}}},\mathbf{s_0}\neq \mathbf{i}}K_{1}\left(\frac{x-X_{\mathbf{i}}}{T}\right)K_{2}\left(h_{\mathbf{n},\mathbf{s_0}}^{-1}\left \|\frac{\mathbf{s_0}-\mathbf{i}}{\mathbf{n}}\right\|\right)}.
\end{equation*}
Let  us set the following sequences, for all $\mathbf{n}\in \mathbb{N}^{*N}$
$$v_{\mathbf{n}}=\left(\frac{k_\mathbf{n}}{k^{'}_{\mathbf{n}}}\right)^{1/d}+\left(\frac{\log(\hat{\mathbf{n}})}{k_\mathbf{n}}\right)^{1/2},$$
and  for all  $\beta \in ]0,1[$ and $x \in D $
\begin{equation}
D_{\mathbf{n}}^{-}(\beta,x)=\left(\frac{k_\mathbf{n}}{cf(x)k^{'}_{\mathbf{n}}}\right)^{1/d}\beta^{1/2d}, \qquad
D_{\mathbf{n}}^{+}(\beta,x)=\left(\frac{k_\mathbf{n}}{cf(x)k^{'}_{\mathbf{n}}}\right)^{1/d}\beta^{-1/2d},
\label{d-+}
\end{equation}
where $c$ is the volume of the unit sphere in $\mathbb{R}^{d}$. It is clear that
$$\forall  \, \mathbf{n}\in \mathbb{N}^{*N}, \forall x \in D \qquad D_{\mathbf{n}}^{-}(\beta,x)\leq D_{\mathbf{n}}^{+}(\beta,x).$$
\begin{lemma}\label{lemtech1}
	If the following conditions are verified:
	\begin{itemize}
		\item[$(L_{1})$]$\displaystyle  \mathbb{I}_{\left\{D_{\mathbf{n}}^{-} (\beta,x)\leq H_{\mathbf{n},x}\leq D_{\mathbf{n}}^{+}(\beta ,x),\, \forall x \in D \right\}} \;\longrightarrow  1 \quad a.co.$
		\item[$(L_{2})$] $\displaystyle \sup_{x \in D}\left \vert \frac{\sum_{\mathbf{i}\in \mathcal{I}_{\mathbf{n}},\mathbf{s_0}\neq \mathbf{i}}K_{1}\left(\frac{x-X_{\mathbf{i}}}{D_{\mathbf{n}}^{-}(\beta,x)}\right) K_{2}\left(h_{\mathbf{n},\mathbf{s_0}}^{-1}\left \|\frac{\mathbf{s_0}-\mathbf{i}}{\mathbf{n}}\right\|\right)}{\sum_{\mathbf{i}\in \mathcal{I}_{\mathbf{n}},\mathbf{s_0}\neq \mathbf{i}}K_{1}\left(\frac{x-X_{\mathbf{i}}}{D_{\mathbf{n}}^{+}(\beta,x)}\right)K_{2}\left(h_{\mathbf{n},\mathbf{s_0}}^{-1}\left \|\frac{\mathbf{s_0}-\mathbf{i}}{\mathbf{n}}\right\|\right)}-\beta\right\vert  \longrightarrow 0\quad a.co.$
		\item[$(L_{3})$]
		$\displaystyle \sup_{x \in D}\left\vert  c_{\mathbf{n}}\left(D_{\mathbf{n}}^{-}(\beta,x)\right)-r(x)\right \vert\longrightarrow 0\quad a.co.$, $\displaystyle\sup_{x\in D} \left \vert c_{\mathbf{n}}\left(D_{\mathbf{n}}^{+}(\beta,x)\right)-r(x) \right\vert \longrightarrow 0\quad a.co$,
	\end{itemize}
	then we have $\displaystyle \sup_{x \in D} \rvert c_{\mathbf{n}}\left(H_{\mathbf{n},x}\right)- r(x)\rvert \longrightarrow 0\qquad a.co.$
\end{lemma}
\begin{lemma}\label{lemtech}
	Under the following conditions:
	\begin{itemize}
		\item[$(L_{1})$]$\displaystyle  \mathbb{I}_{\left\{D_{\mathbf{n}}^{-} (\beta,x)\leq H_{\mathbf{n},x}\leq D_{\mathbf{n}}^{+}(\beta ,x),\, \forall x \in D \right\}} \;\longrightarrow  1 \quad a.co.$
		\item[$(L_{2}^{'})$] $\displaystyle \sup_{x \in D}\left \vert \frac{\sum_{\mathbf{i}\in \mathcal{I}_{\mathbf{n}},\mathbf{s_0}\neq \mathbf{i}}K_{1}\left(\frac{x-X_{\mathbf{i}}}{D_{\mathbf{n}}^{-}(\beta,x)}\right) K_{2}\left(h_{\mathbf{n},\mathbf{s_0}}^{-1}\left \|\frac{\mathbf{s_0}-\mathbf{i}}{\mathbf{n}}\right\|\right)}{\sum_{\mathbf{i}\in \mathcal{I}_{\mathbf{n}},\mathbf{s_0}\neq \mathbf{i}}K_{1}\left(\frac{x-X_{\mathbf{i}}}{D_{\mathbf{n}}^{+}(\beta,x)}\right)K_{2}\left(h_{\mathbf{n},\mathbf{s_0}}^{-1}\left \|\frac{\mathbf{s_0}-\mathbf{i}}{\mathbf{n}}\right\|\right)}-\beta\right\vert =\mathcal{O}(v_{\mathbf{n}}) \qquad a.co.$
		\item[$(L_{3}^{'})$]
		$\displaystyle \sup_{x \in D}\left\vert  c_{\mathbf{n}}\left(D_{\mathbf{n}}^{-}(\beta,x)\right)-r(x)\right \vert=\mathcal{O}(v_{\mathbf{n}}) \qquad a.co$,\\  $\displaystyle\sup_{x\in D} \left \vert c_{\mathbf{n}}\left(D_{\mathbf{n}}^{+}(\beta,x)\right)-r(x) \right\vert =\mathcal{O}(v_{\mathbf{n}})\qquad a.co,$
	\end{itemize}
	we have,  $\displaystyle \sup_{x \in D}\left \vert c_{\mathbf{n}}\left(H_{\mathbf{n},x}\right)- r(x)\right\vert=\mathcal{O}(v_{\mathbf{n}})\qquad a.co.$
\end{lemma}
\noindent The proof of Lemma~\ref{lemtech} is similar as in \cite{col80} and is therefore omitted. Lemma~\ref{lemtech1} is a particular case of the proof of Lemma~\ref{lemtech} when we take  $v_{\mathbf{n}}=1$ and $C=1$.
\subsection*{Proofs of Lemma~\ref{th1_reg} and Lemma~\ref{th2_reg}}
\noindent  Since the proof of Lemma \ref{th1_reg} is based on the result of Lemma \ref{lemtech1}, it is  sufficient  to check conditions $(L_{1})$, $(L_{2})$ and $(L_{3})$. For the proof of Lemma~\ref{th2_reg}, it suffices to check conditions $(L_{2}')$ and  $(L_{3}')$.

\noindent To check the condition $ (L_1 )$, we need the following two lemmas.

\begin{lemma}\label{l12}
	(\cite{ibra} or \cite{deo})
	\begin{itemize}
		\item [i)] We assume that the mixing condition (\ref{cddep}) is satisfied. We denote by $\mathcal{L}_{r}\left(\mathcal{F}\right)$ the class of $\mathcal{F}-$mesurable random variables $X$ satisfying $$\|X\|_{r}:=\left(E\left(\lvert X\rvert^{r}\right)\right)^{1/r}< \infty.$$ Let $X\in \mathcal{L}_{r}\left(\mathcal{B}(E)\right)$, $Y\in \mathcal{L}_{s}\left(\mathcal{B}(E')\right)$ and 
		$ 1\leq r,s,t\leq \infty$ such that
		$\frac{1}{r}+\frac{1}{s}+\frac{1}{t}=1$, then
		\begin{equation}
		\lvert \mathrm{Cov}(X,Y)\rvert \leq  \|X\|_{r}\|Y\|_{s}\left\{\psi\left(\mathrm{Card(E),Card(E^{'})}\right)\varphi\left(\mathrm{dist}(E,E^{'})\right)\right\}^{1/t}.
		\label{Ibli1}
		\end{equation}
		\item [ii)] For random variables $X,Y$ bounded with probability $1$, we have
		\begin{equation}
		 \lvert \mathrm{Cov}(X,Y)\rvert \leq C\psi\left(\mathrm{Card(E),Card(E^{'})}\right)\varphi\left(\mathrm{dist}(E,E^{'})\right).
		\label{Ibli2}
		\end{equation}
	\end{itemize}
\end{lemma}

\begin{lemma}\label{l1}
	Under assumptions of  Theorem \ref{th1}, we have
	$$S_{\mathbf{n}}+R_{\mathbf{n}}=\mathcal{O}\left(k^{'}_{\mathbf{n}}\delta_{\mathbf{n}}\right),$$
	where
	$$ S_{\mathbf{n}}= \sum_{\mathbf{i} \in \mathcal{V}_{\mathbf{s_0}}}\mathrm{Var}\left(\Lambda_{\mathbf{i}}\right) \; \mathrm{and}
	\quad R_{\mathbf{n}}=\sum_{\mathbf{i} \in \mathcal{V}_{\mathbf{s_0}}}\sum_{\underset{\mathbf{j}\neq \mathbf{i}}{\mathbf{j} \in \mathcal{V}_{\mathbf{s_0}}}}\left|\mathrm{Cov}\left(\Lambda_{\mathbf{i}},\Lambda_{\mathbf{j}}\right)\right|,$$
	$$ \Lambda_{\mathbf{i}}=\mathbb{I}_{B(x,D_{\mathbf{n}})}(X_{\mathbf{i}}), \quad \mathbf{i}\in \mathcal{I}_{\mathbf{n}}, \quad \delta_{\mathbf{n}}=\mathbb{P}\left(\|X-x\|<D_{\mathbf{n}}\right),\quad  D_{\mathbf{n}}^{d}=\mathcal{O}\left(\frac{k_\mathbf{n}}{k^{'}_{\mathbf{n}}}\right),$$
	$B(x,\varepsilon)$ denotes  the closed ball of\,  $\mathbb{R}^{d}$ with center $x$ and  radius $\varepsilon$.
\end{lemma}
\subsection*{Proof of Lemma \ref{l1}}
\noindent Let  $\delta_{\mathbf{n},\mathbf{i}}=\mathbb{P}\left(\|X_{\mathbf{i}}-x\|< D_{\mathbf{n}}\right)$, we can deduce that
$$ S_{\mathbf{n}}=\sum_{\mathbf{i} \in \mathcal{V}_{\mathbf{s_0}}}Var\left(\Lambda_{\mathbf{i}}\right)=\sum_{\mathbf{i} \in \mathcal{V}_{\mathbf{s_0}}}\delta_{\mathbf{n},\mathbf{i}}(1-\delta_{\mathbf{n},\mathbf{i}})=\mathcal{O}\left(k^{'}_{\mathbf{n}}\delta_{\mathbf{n}}\right),$$
by  the following results. \\ Firstly, under the Lipschitz condition of $f$ (assumption (H1)), we have
\begin{eqnarray}
\delta_{\mathbf{n}}&=&\mathbb{P}\left(\|X-x\|<  D_{\mathbf{n}}\right) \nonumber\\
&=&f(x)\int _{B(x,D_{\mathbf{n}})}du+\int _{B(x,D_{\mathbf{n}})} (f(u)-f(x))du\nonumber\\
&=& c f(x) D_{\mathbf{n}}^{d}+ \mathcal{O}\left(D_{\mathbf{n}}^{d+1} \right).
\label{eq3}
\end{eqnarray}
Secondly
\begin{eqnarray}
\delta_{\mathbf{n},\mathbf{i}}-\delta_{\mathbf{n}}&=&\int_{B(x,D_{\mathbf{n}})}\left(f_{\mathbf{i}}(u)-f(u)\right)(u)du \nonumber\\
&=&\sup_u\left|f_{\mathbf{i}}(u)-f(u)\right|\mathcal{O}\left(\frac{k_{\mathbf{n}}}{k^{'}_{\mathbf{n}}}\right).
\end{eqnarray}
Thus, the local stationarity assumption ($H_8$) implies  
\begin{equation}
\sum_{\mathbf{i}\in \mathcal{V}_{\mathbf{s_0}}}\left(\delta_{\mathbf{n},\mathbf{i}}-\delta_{\mathbf{n}}\right)=o(k_{\mathbf{n}}).
\label{deltai}
\end{equation}
Now for $\displaystyle R_{\mathbf{n}}$, it should be noted that by (H5) and for each  $\mathbf{j}\neq \mathbf{i}$
\begin{eqnarray}
\left \vert \mathrm{Cov}\left(\Lambda_{\mathbf{i}},\Lambda_{\mathbf{j}}\right)\right \vert &=&\left \vert \mathbb{P}\left(\|X_{\mathbf{i}}-x\|<  D_{\mathbf{n}},\|X_{\mathbf{j}}-x\|<  D_{\mathbf{n}}\right)\right.\nonumber\\
&\qquad &\qquad -  \left. \mathbb{P}\left(\Vert X_{\mathbf{i}}-x \Vert< D_{\mathbf{n}}\right) \mathbb{P}\left( \Vert X_{\mathbf{j}}-x \Vert <  D_{\mathbf{n}}\right)\right\vert \nonumber\\
&\leq&\int_{B\left(x,D_{\mathbf{n}}\right)\times B\left(x,D_{\mathbf{n}}\right)} \left \vert f_{X_{\mathbf{i}}X_{\mathbf{j}}}(u,v)-f_{\mathbf{i}}(u)f_{\mathbf{j}}(v)\right \vert dudv \nonumber\\
&\leq & C D_{\mathbf{n}}^{2d} \leq C \delta_\mathbf{n}^2,
\label{eq31}
\end{eqnarray}
since by (\ref{eq3})
$$\frac{D_{\mathbf{n}}^{d}}{\delta_{\mathbf{n}}}\to \frac{1}{cf(x)},\qquad \mathrm{as}\quad \mathbf{n}\to \infty.$$
Using  Lemma \ref{l12} and (\ref{eq3}),  we can write for $\displaystyle r=s= 4 $
\begin{eqnarray}
\left|\mathrm{Cov}\left(\Lambda_{\mathbf{i}},\Lambda_{\mathbf{j}}\right)\right| &\leq&C\left[E\left(\Lambda_{\mathbf{i}}^{4}\right)E\left(\Lambda_{\mathbf{j}}^{4}\right)\right]^{1/4}\left(\psi(1,1)\varphi\left(\|\mathbf{i}-\mathbf{j}\|\right)\right)^{1/2}\nonumber\\
&\leq& C\delta_{\mathbf{n}}^{1/2}\varphi\left(\|\mathbf{i}-\mathbf{j}\|\right)^{1/2}.
\label{eq32}
\end{eqnarray}
Let $q_{\mathbf{n}}$ be a sequence of real numbers  defined by $\displaystyle q_{\mathbf{n}}^{N}=\mathcal{O}\left(\frac{k^{'}_{\mathbf{n}}}{k_\mathbf{n}}\right)$. Using the later, we  define $\displaystyle S=\{ \mathbf{i},\mathbf{j}\in
\mathcal{V}_{\mathbf{s_0}},\; 0<\|\mathbf{i}-\mathbf{j}\|\leq q_{\mathbf{n}}\}$ and $S^{c}$ its complementary in
$\mathcal{V}_{\mathbf{s_0}}$, and rewrite
 $$ R_{\mathbf{n}}=\sum_{\mathbf{i},\mathbf{j}\in S}\left |\mathrm{Cov}\left(\Lambda_{\mathbf{i}},\Lambda_{\mathbf{j}}\right)\right|+\sum_{\mathbf{i},\mathbf{j}\in S^{c}}\left |\mathrm{Cov}\left(\Lambda_{\mathbf{i}},\Lambda_{\mathbf{j}}\right)\right|= R_{\mathbf{n}}^{(1)}+R_{\mathbf{n}}^{(2)}.$$
\noindent Firstly, according to the definitions of $q_{\mathbf{n}}$ and $S$, and equation (\ref{eq31}), we have
\begin{eqnarray*}
	R_{\mathbf{n}}^{(1)}&\leq &\sum_{\mathbf{i},\mathbf{j}\in S}C\delta_{\mathbf{n}}^{2}
	\leq C\delta_{\mathbf{n}}^{2}k^{'}_{\mathbf{n}}q_{\mathbf{n}}^{N}
	=\mathcal{O}\left(k^{'}_{\mathbf{n}}\delta_{\mathbf{n}}\right),
\end{eqnarray*}
since $\displaystyle\delta_{\mathbf{n}}=\mathcal{O}(q_{\mathbf{n}}^{-N})$ by (\ref{eq3}).\\
\noindent Secondly, by (\ref{cd3}) and  (\ref{eq32}), we get
\begin{eqnarray*}
	R_{\mathbf{n}}^{(2)}&\leq &C\delta_{\mathbf{n}}^{1/2}\sum_{\mathbf{i},\mathbf{j}\in S^{c}}\varphi\left(\|\mathbf{i}-\mathbf{j}\|\right)^{1/2}
	=C\delta_{\mathbf{n}}^{1/2}k^{'}_{\mathbf{n}}\sum_{\mathbf{i}\in S^{c}}\varphi\left(\|\mathbf{i}\|\right)^{1/2}\\
	&=&C\delta_{\mathbf{n}}k^{'}_{\mathbf{n}}\delta_{\mathbf{n}}^{-1/2}\sum_{\mathbf{i}\in S^{c}}\varphi\left(\|\mathbf{i}\|\right)^{1/2}
	\leq C\delta_{\mathbf{n}}k^{'}_{\mathbf{n}}\left(\frac{k_\mathbf{n}}{k^{'}_{\mathbf{n}}}\right)^{-1/2}\sum_{\mathbf{i}\in S^{c}}\varphi\left(\|\mathbf{i}\|\right)^{1/2}\\
	&\leq & C\delta_{\mathbf{n}}k^{'}_{\mathbf{n}}\sum_{\mathbf{i}\in S^{c}}\|\mathbf{i}\|^{(N-\theta)/2}
 =\mathcal{O}\left(\delta_{\mathbf{n}}k^{'}_{\mathbf{n}}\right),
\end{eqnarray*}
because under  assumptions $(H6)$ and $(H7)$, we have $\theta>(1+\frac{2\gamma}{\gamma-\tilde{\gamma}})N$, thus
\begin{equation*}
\sum_{\mathbf{i}\in S^{c}}\|\mathbf{i}\|^{(N-\theta)/2}\leq k^{'}_{\mathbf{n}}q_{\mathbf{n}}^{^{(N-\theta)/2}} = o(1).
\end{equation*}
Finally, the result follows:  $$R_{\mathbf{n}}=\mathcal{O}\left(k^{'}_{\mathbf{n}}\delta_{\mathbf{n}}\right)
\; \mathrm{and}\;  S_{\mathbf{n}}+R_{\mathbf{n}}=\mathcal{O}\left(k^{'}_{\mathbf{n}}\delta_{\mathbf{n}}\right).$$

\subsection*{Verification of $(L_{1})$}
\noindent Let $\varepsilon_\mathbf{n}=\frac{1}{2} \varepsilon_0\left(k_\mathbf{n}/k^{'}_{\mathbf{n}}\right)^{1/d}$ with $\varepsilon_{0}>0$ and let $N_{\varepsilon_\mathbf{n}}=\mathcal{O}(\varepsilon_{\mathbf{n}}^{-d})$ be a positive integer. Since $D$ is compact, one can cover it by   $N_{\varepsilon_\mathbf{n}}$  closed balls in $\mathbb{R}^d$ of centers $x_i \in D, \; i=1,\ldots, N_{\varepsilon_\mathbf{n}}$ and radius $\varepsilon_\mathbf{n}$. Let us show that
$$\mathbb{I}_{\left\{D_{\mathbf{n}}^{-}(\beta, x) \leq H_{\mathbf{n},x}\leq D_{\mathbf{n}}^{+}(\beta, x),\, \forall x \in D\right\}}\longrightarrow 1 \qquad a.co,$$
which can be written as, $\forall \; \eta > 0$,
$$\sum_{\mathbf{n}\in\mathbb{N}^{*N}}\mathbb{P}(\mid \mathbb{I}_{\left\{D_{\mathbf{n}}^{-}(\beta, x) \leq H_{\mathbf{n},x}\leq D_{\mathbf{n}}^{+}(\beta, x),\, \forall x \in D\right \}} -1\mid > \eta) < \infty. $$
We have
\begin{align}
\,& \mathbb{P}(\mid \mathbb{I}_{\left\{D_{\mathbf{n}}^{-}(\beta, x) \leq H_{\mathbf{n},x}\leq D_{\mathbf{n}}^{+}(\beta, x),\, \forall x \in D\right \}} -1 \mid > \eta) \nonumber \\
&\leq  \mathbb{P}\left(\sup _{x \in D}\left(H_{\mathbf{n},x} - D_{\mathbf{n}}^{-}(\beta,x)\right) <0 \right) + \mathbb{P}\left(\inf_{x \in D} \left(H_{\mathbf{n},x} -D_{\mathbf{n}}^{+}(\beta,x) \right)> 0\right)\nonumber \\
&\leq  \mathbb{P}\left(\max _{1\leq i \leq N_{\varepsilon_\mathbf{n}}} \left(H_{\mathbf{n},x_i} - D_{\mathbf{n}}^{-}(\beta,x_i)\right) <2\varepsilon_\mathbf{n} \right) + \mathbb{P}\left(\min_{1\leq i \leq N_{\varepsilon_\mathbf{n}}}\left( H_{\mathbf{n},x_i} -D_{\mathbf{n}}^{+}(\beta,x_i) \right)> -2\varepsilon_\mathbf{n}\right) \nonumber \\
&\leq N_{\varepsilon_\mathbf{n}} \max_{1\leq i \leq N_{\varepsilon_\mathbf{n}}} \mathbb{P}\left( H_{\mathbf{n},x_i} <  D_{\mathbf{n}}^{-}(\beta,x_i)+ 2\varepsilon_\mathbf{n} \right)  \nonumber \\
& \qquad \qquad \qquad + N_{\varepsilon_\mathbf{n}} \max_{1\leq i \leq N_{\varepsilon_\mathbf{n}}} \mathbb{P}\left( H_{\mathbf{n},x_i}  > D_{\mathbf{n}}^{+}(\beta,x_i) -2\varepsilon_\mathbf{n}\right).
\label{eq1}
\end{align}
Let us evaluate the first term in the right-hand side of (\ref{eq1}), without ambiguity we ignore the $i$ index  in $x_i$. As justified in the following
\begin{eqnarray}
\mathbb{P}\left(H_{\mathbf{n},x} <D_{\mathbf{n}}^{-}(\beta,x)+2\varepsilon_\mathbf{n} \right)
& \leq &\mathbb{P}\left(\sum_{\mathbf{i} \in \mathcal{V}_{\mathbf{s_0}}}\mathbb{I}_{B(x,D_{\mathbf{n}}^{-}(\beta,x)+2\varepsilon_\mathbf{n})}(X_{\mathbf{i}}) >k_\mathbf{n}\right) \label{eqp1} \\
&\leq & \mathbb{P}\left(\sum_{\mathbf{i} \in \mathcal{V}_{\mathbf{s_0}}}\xi_{\mathbf{i}} >k_\mathbf{n}-k^{'}_{\mathbf{n}}\delta_{\mathbf{n}}\right)
\label{eqp2} \\
&\leq &\mathbb{P}\left(\sum_{\mathbf{i} \in \mathcal{V}_{\mathbf{s_0}}}\xi_{\mathbf{i}} >Ck_\mathbf{n}(1-\beta^{1/2})\right)
:= P_{1,\mathbf{n}}, \label{eqp3}
\end{eqnarray}
where $\displaystyle \xi_{\mathbf{i}}=\Lambda_\mathbf{i}-\delta_{\mathbf{n},\mathbf{i}}$ is centered and  $\Lambda_\mathbf{i}$ is defined in Lemma \ref{l1} when we replace  $D_{\mathbf{n}}$ by $D_{\mathbf{n}}^{-}+2\varepsilon_\mathbf{n}$. From (\ref{eqp1}), we get (\ref{eqp2}) by (\ref{deltai}) while result (\ref{eqp2}) permits to get (\ref{eqp3}) by the help of the following.

\noindent Actually, according to the definition of $D_{\mathbf{n}}^{-}$  in (\ref{d-+}) and  replacing $D_{\mathbf{n}}$ by $D_{\mathbf{n}}^{-}+2\varepsilon_\mathbf{n}$ in (\ref{eq3}), we get
\begin{eqnarray}
k^{'}_{\mathbf{n}}\delta_{\mathbf{n}}-k_\mathbf{n}\left(\varepsilon_0 (cf(x))^{1/d} + \beta^{1/2d}\right)^{d}=o(k_\mathbf{n}),
\label{eq4}
\end{eqnarray}
therefore, for all $\varepsilon_{1}>0$,
\begin{eqnarray*}
	k_\mathbf{n} - k^{'}_{\mathbf{n}}\delta_{\mathbf{n}}>
	k_\mathbf{n}\left( 1-\left(\varepsilon_0 (c f(x))^{1/d} + \beta^{1/2d}\right)^{d} - \varepsilon_1\right).
\end{eqnarray*}
Then, for $\varepsilon_1$ and $\varepsilon_0$ very small such that  $ 1-\left(\varepsilon_0 (cf(x))^{1/d} + \beta^{1/2d}\right)^{d} - \varepsilon_1>0$, we can  find some constant  $C>0$ such that
\begin{eqnarray}
k_\mathbf{n}-k^{'}_{\mathbf{n}}\delta_{\mathbf{n}}>C k_\mathbf{n}(1-\beta^{1/2}).
\label{eq5}
\end{eqnarray}
For the second term in the right-hand side of (\ref{eq1}),
\begin{eqnarray}
\mathbb{P}\left(H_{\mathbf{n},x}>D_{\mathbf{n}}^{+}(\beta ,x)-2\varepsilon_\mathbf{n} \right)&\leq &\mathbb{P}\left(\sum_{\mathbf{i} \in \mathcal{V}_{\mathbf{s_0}}}\mathbb{I}_{B(x,D_{\mathbf{n}}^{+}(\beta,x)-2\varepsilon_\mathbf{n})}(X_{\mathbf{i}})<k_\mathbf{n}\right)
\label{eqp21} \\
&\leq& \mathbb{P}\left(\sum_{\mathbf{i} \in \mathcal{V}_{\mathbf{s_0}}}\Delta_{\mathbf{i}} >k^{'}_{\mathbf{n}}\delta_{\mathbf{n}}-k_\mathbf{n}\right)
\label{eqp22} \\
&\leq & \mathbb{P}\left(\sum_{\mathbf{i} \in \mathcal{V}_{\mathbf{s_0}}}\Delta_{\mathbf{i}} >C k_\mathbf{n}\left(\beta^{-1/2}-1\right)\right)
:=P_{2,\mathbf{n}}, \label{eqp23}
\end{eqnarray}
where $\displaystyle \Delta_{\mathbf{i}}=\delta_{\mathbf{n},\mathbf{i}}-\Lambda_{\mathbf{i}}$ is centered and $\Lambda_\mathbf{i}$ is defined in Lemma \ref{l1} replacing  $D_{\mathbf{n}}$ by $D_{\mathbf{n}}^{+} - 2\varepsilon_\mathbf{n}$. Result  (\ref{eqp22}) is obtained by  (\ref{deltai}) while that of (\ref{eqp23}) is obtained by replacing $D_{\mathbf{n}}$ in (\ref{eq3})  by $D_{\mathbf{n}}^{+}-2\varepsilon_\mathbf{n}$. Then, we get 
\begin{equation}
k^{'}_{\mathbf{n}}\delta_\mathbf{n}- k_\mathbf{n} \left(\beta^{-1/2d}-\varepsilon_0(cf(x))^{1/d}\right)^{d}=o(k_\mathbf{n}).
\label{eq6}
\end{equation}
Thus for all $\varepsilon_{2}>0$, it is easy to see that
$$k^{'}_{\mathbf{n}}\delta_{\mathbf{n}}-k_\mathbf{n}>k_\mathbf{n}\left(\left(\beta^{-1/2d}-\varepsilon_0 (cf(x))^{1/d}\right)^{d}-1-\varepsilon_{2}\right),$$
so  for $\varepsilon_{2}$ and $\varepsilon_0$ small  enough such that  $\left(\left(\beta^{-1/2d}-\varepsilon_0 (cf(x))^{1/d}\right)^{d}-1-\varepsilon_{2}\right)>0$, there exists  $C>0$ such that
\begin{equation}
k^{'}_{\mathbf{n}}\delta_{\mathbf{n}}-k_\mathbf{n}>C k_\mathbf{n}\left(\beta^{-1/2}-1\right).
\end{equation}
Now, it suffices to prove that $$ \sum_{\mathbf{n}\in\mathbb{N}^{*N}}N_{\varepsilon_{\mathbf{n}}}P_{1,\mathbf{n}}<\infty \quad \mathrm{and}\quad\sum_{\mathbf{n}\in\mathbb{N}^{*N}}N_{\varepsilon_{\mathbf{n}}}P_{2,\mathbf{n}}<\infty.$$
\subsubsection*{Let us consider $P_{1,\mathbf{n}}$}
This proof is based on the classical spatial block decomposition of the sum on
$\xi_{\mathbf{i}}$ in $\mathcal{V}_{\mathbf{s_0}}$ similarly to  \cite{tran90}.
Let $\mathcal{G}_\mathbf{n}\subset \mathcal{I}_{\mathbf{n}}$ be the smallest rectangular grid of center $\mathbf{s_0}$ containing  $\mathcal{V}_{\mathbf{s_0}}$. Without loss of generality, we assume that $\mathcal{G}_\mathbf{n}$ is defined via some $\mathbf{k}=(k_1,\ldots,k_N)$ where $1\leq k_j\leq n_j, j=1,\ldots,N$. However, by construction $\mathcal{G}_\mathbf{n}$ is of cardinal  $\hat{\mathbf{k}}=k_1\times\cdots \times k_N$  satisfying $k^{'}_{\mathbf{n}}=\mathcal{O}(\hat{\mathbf{k}})$. In addition, we assume that $k_{l}=2bq_{l} \; ,\; l=1,\ldots,N $, where $q_{l}$ and $b$ are positive integers. Then  the decomposition can be presented as follows
\begin{eqnarray*}
	U(1,\mathbf{k},\mathbf{j})&=&\sum_{\underset{k=1,\ldots,N.}{i_{l}=2j_{l}b+1,}}^{(2j_{l}+1)b}\xi_{\mathbf{i}} \\
	U(2,\mathbf{k},\mathbf{j})&=&\sum_{\underset{l=1,\ldots,N-1.}{i_{l}=2j_{l}b+1,}}^{(2j_{l}+1)b} \; \; \sum_{i_{N}=(2j_{N}+1)b+1,}^{2(j_{N}+1)b}\xi_{\mathbf{i}} \\
	U(3,\mathbf{k},\mathbf{j})&=&\sum_{\underset{l=1,\ldots,N-2.}{i_{l}=2j_{l}b+1,}}^{(2j_{l}+1)b}\; \; \sum_{i_{N-1}=(2j_{N-1}+1)b+1,}^{2(j_{N-1}+1)b} \; \; \sum_{i_{N}=2j_{N}b+1,}^{(2j_{N}+1)b}\xi_{\mathbf{i}}
	\\
	U(4,\mathbf{k},\mathbf{j})&=&\sum_{\underset{l=1,\ldots,N-2.}{i_{l}=2j_{l}b+1,}}^{(2j_{l}+1)b}\; \; \sum_{i_{N-1}=(2j_{N-1}+1)b+1,}^{2(j_{N-1}+1)b} \; \; \sum_{i_{N}=(2j_{N}+1)b+1,}^{2(j_{N}+1)b}\xi_{\mathbf{i}}\\
	...
\end{eqnarray*}
Note that
$$ U(2^{N-1},\mathbf{k},\mathbf{j})=\sum_{\underset{l=1,\ldots,N-1.}{i_{l}=(2j_{l}+1)b+1,}}^{2(j_{l}+1)b} \; \; \sum_{i_{N}=2j_{N}b+1,}^{(2j_{N}+1)b}\xi_{\mathbf{i}}$$
and that
$$U(2^{N},\mathbf{k},\mathbf{j})=\sum_{\underset{l=1,\ldots,N.}{i_{l}=(2j_{l}+1)b+1,}}^{2(j_{l}+1)b}\xi_{\mathbf{i}}.$$
For each integer  $1\leq l \leq 2^{N}$, let
$$ T(\mathbf{k},l)=\sum_{\underset{l=1,\ldots,N}{j_{l}=0}}^{q_{l}-1} U(l,\mathbf{k},\mathbf{j}).$$
Therefore, we have
\begin{eqnarray}
\sum_{\mathbf{i} \in \mathcal{V}_{\mathbf{s_0}}}\xi_{\mathbf{i}}=\sum_{l=1}^{2^{N}}T(\mathbf{k},l).
\label{eq7}
\end{eqnarray}
It follows that
\begin{eqnarray*}
	P_{1,\mathbf{n}}&=&\mathbb{P}\left (\sum_{l=1}^{2^{N}}T(\mathbf{k},l)>Ck_\mathbf{n}(1-\sqrt{\beta})\right)\leq 
	2^{N}\mathbb {P}\left (\mid T(\mathbf{k},1)\mid>\frac{Ck_\mathbf{n}(1-\sqrt{\beta})}{2^{N}}\right).
\end{eqnarray*}
We enumerate in an arbitrary manner the  $\hat{\mathbf{q}}=q_{1}\times \ldots \times q_{N}$ terms $\displaystyle U(1,\mathbf{k},\mathbf{j})$ of the sum $T(\mathbf{k},1)$ and denote them  $\displaystyle W_{1},\ldots,W_{\hat{\mathbf{q}}}$.
Notice that, $\displaystyle U(1,\mathbf{k},\mathbf{j})$ is measurable with respect to the field generated by the  $Z_{\mathbf{i}}$ with  $\mathbf{i}\in \mathbf{I}(\mathbf{k},\mathbf{j})=\{ \mathbf{i}\in \mathcal{G}_{\mathbf{n}} \; | \;2j_{l}b+1\leq i_{l}\leq(2j_{l}+1)b ,\;l=1,\ldots,N\}$, the set  $\mathbf{I}(\mathbf{k},\mathbf{j})$ contains  $b^{N}$ sites and $\mathrm{dist}(\mathbf{I}(\mathbf{k},\mathbf{j}),\mathbf{I}(\mathbf{k},\mathbf{j}^{'}))>b$. In addition, we have
$\mid W_{l}\mid\leq b^{N}.$

\noindent According to Lemma 4.5 of \cite{carbon1997}, one can find a sequence of independent random variables $\displaystyle W_{1}^{*},\ldots,W_{\hat{\mathbf{q}}}^{*}$ where $W_{l}$ has the same distribution as  $W_{l}^{*}$ and:
\begin{eqnarray*}
	\sum_{l=1}^{\hat{\mathbf{q}}}\Ex(\mid W_{l}-W_{l}^{*}\mid)&\leq&4\hat{\mathbf{q}}b^{N}\psi((\hat{\mathbf{q}}-1)b^{N},b^{N})\varphi(b).
\end{eqnarray*}
Then, we can write
\begin{eqnarray*}
	P_{1,\mathbf{n}}&\leq &2^{N}\mathbb {P}\left (\mid T(\mathbf{n},1)\mid>\frac{Ck_\mathbf{n}(1-\sqrt{\beta})}{2^{N}}\right)\\
	&\leq & 2^{N}\mathbb {P}\left (\mid \sum_{l=1}^{\hat{\mathbf{q}}}W_{l}\mid>\frac{Ck_\mathbf{n}(1-\sqrt{\beta})}{2^{N}}\right)\\
	&\leq&2^{N}\mathbb {P}\left ( \sum_{l=1}^{\hat{\mathbf{q}}}\mid W_{l}-W_{l}^{*}\mid>\frac{Ck_\mathbf{n}(1-\sqrt{\beta})}{2^{N+1}}\right)\\
	&\, & \qquad +2^{N}\mathbb {P}\left ( \sum_{l=1}^{\hat{\mathbf{q}}}\mid W_{l}^{*}\mid>\frac{Ck_\mathbf{n}(1-\sqrt{\beta})}{2^{N+1}}\right).
\end{eqnarray*}
Let
$P_{11,\mathbf{n}}=\mathbb {P}\left ( \sum_{l=1}^{\hat{\mathbf{q}}}\mid W_{l}-W_{l}^{*}\mid>\frac{Ck_\mathbf{n}(1-\sqrt{\beta})}{2^{N+1}}\right)$
and
$P_{12,\mathbf{n}}=\mathbb{P}\left ( \sum_{l=1}^{\hat{\mathbf{q}}}\mid W_{l}^{*}\mid>\frac{Ck_\mathbf{n}(1-\sqrt{\beta})}{2^{N+1}}\right).$ \\
It suffices to show that  $\displaystyle \sum_{\mathbf{n}\in\mathbb{N}^{*N}}P_{11,\mathbf{n}}<\infty$ and $\displaystyle \sum_{\mathbf{n}\in\mathbb{N}^{*N}}P_{12,\mathbf{n}}<\infty$.

\subsubsection*{Let us consider first $P_{11,\mathbf{n}}$}
\noindent Using Markov's inequality, we get
\begin{eqnarray*}
	P_{11,\mathbf{n}}&=&\mathbb {P}\left ( \sum_{l=1}^{\hat{\mathbf{q}}}\mid W_{l}-W_{l}^{*}\mid>\frac{C k_\mathbf{n}(1-\sqrt{\beta})}{2^{N+1}}\right)\\
	&\leq&\frac{2^{N+3}}{C k_\mathbf{n}(1-\sqrt{\beta})}\hat{\mathbf{q}}b^{N}\psi((\hat{\mathbf{q}}-1)b^{N},b^{N})\varphi(b) \\
	&\leq& \frac{C}{k_\mathbf{n}(1-\sqrt{\beta})}k^{'}_{\mathbf{n}}\psi((\hat{\mathbf{q}}-1)b^{N},b^{N})\varphi(b),
\end{eqnarray*}
because $\hat{\mathbf{k}}=2^N\hat{\mathbf{q}}b^{N}$ by definition and $k^{'}_{\mathbf{n}}=\mathcal{O}(\hat{\mathbf{k}})$\\ 
Let us consider that
\begin{equation}
b^{N}=\mathcal{O}\left(\hat{\mathbf{n}}^{2(1-s(1-\tilde{\gamma}))/a}\right),
\label{bchoi}
\end{equation}
where $a=2+(2+s(2-\tilde{\gamma}))d+s(4+2\tilde{\beta}+2\gamma-3\tilde{\gamma})$.\\
Under the assumption on the function $\psi(n,m)$, we distinguish the following two cases:\medskip
\newline \noindent\textbf{Case 1}
$$\psi(n,m)\leq C \min(n,m)\, \text{with}\,   (1-s(1-\tilde{\gamma}))\theta>N\left\{(2+s(2-\tilde{\gamma}))d+2s(2+\gamma-\tilde{\gamma})\right\},$$
$$
\text{ and }\,  2<s<\frac{1}{1-\tilde{\gamma}}.
$$
\\ \medskip
\noindent In this case, we have
\begin{eqnarray*}
	P_{11,\mathbf{n}}&\leq&C\frac{k^{'}_{\mathbf{n}}}{k_{\mathbf{n}}}b^{N}\varphi(b)
	\leq C \frac{k^{'}_{\mathbf{n}}}{k_{\mathbf{n}}}b^{N-\theta}.
\end{eqnarray*}
Then by using (\ref{bchoi}) and the definition  of $N_{\varepsilon_\mathbf{n}}$, we have
$$N_{\varepsilon_\mathbf{n}}P_{11,\mathbf{n}}\leq C \hat{\mathbf{n}}^{-2\left(1-\frac{3+s(2\tilde{\beta}-1)}{a}\right)}.$$
One can show that $a>2(3+s(2\tilde{\beta}-1))$ and then $\displaystyle \sum_{\mathbf{n}\in\mathbb{N}^{*N}}N_{\varepsilon_\mathbf{n}}P_{11,\mathbf{n}}<\infty$.\\ \medskip
\noindent\textbf{Case 2}
$$\psi(n,m)\leq C(n+m+1)^{\tilde{\beta}}\,\text{with}\, (1-s(1-\tilde{\gamma}))\theta>N\left\{2+(2+s(2-\tilde{\gamma}))d+s(4+2\tilde{\beta}+2\gamma-3\tilde{\gamma})\right\} \text{and }$$
$2<s<\frac{1}{1-\tilde{\gamma}}$. In this case, we have
\begin{eqnarray*}
	P_{11,\mathbf{n}}&\leq&C\frac{k^{'}_{\mathbf{n}}}{k_\mathbf{n}}(k^{'}_{\mathbf{n}}b^{N})^{\tilde{\beta} }\varphi(b)
	\leq C\frac{k_{\mathbf{n}}^{'}}{k_\mathbf{n}} k_{\mathbf{n}}^{'\tilde{\beta}}b^{-\theta}\leq C \hat{\mathbf{n}}^{-(2-\gamma \tilde{\beta})}.
\end{eqnarray*}
Then, it follows that
$\displaystyle \sum_{\mathbf{n}\in\mathbb{N}^{N}}N_{\varepsilon_\mathbf{n}}P_{11,\mathbf{n}}<\infty$ when $\tilde{\beta}<1/\gamma$. 

\subsubsection*{Let us consider $P_{12,\mathbf{n}}$}
\noindent Applying Markov's inequality, we have for $t>0$:
\begin{eqnarray*}
	P_{12,\mathbf{n}}&=&\mathbb {P}\left ( \sum_{l=1}^{\hat{\mathbf{q}}}\mid W_{l}^{*}\mid>\frac{Ck_\mathbf{n}(1-\sqrt{\beta})}{2^{N+1}}\right)\\
	&\leq&\exp\left(-t\frac{Ck_\mathbf{n}(1-\sqrt{\beta})}{2^{N+1}}\right)\Ex\left(\exp\left(t\sum_{l=1}^{\hat{\mathbf{q}}}W_{l}^{*}\right)\right)\\
	&\leq&\exp\left(-t\frac{Ck_\mathbf{n}(1-\sqrt{\beta})}{2^{N+1}}\right)\prod_{l=1}^{\hat{\mathbf{q}}}\Ex\left(\exp\left(tW_{l}^{*}\right)\right),
\end{eqnarray*}
since the variables $W_{1}^{*},\ldots,W_{\hat{\mathbf{q}}}^{*}$ are independent.

\noindent Let $r>0 $, for $\displaystyle t=\frac{r\log(\hat{\mathbf{n}})}{k_\mathbf{n}}$, $l=1,\ldots,\hat{\mathbf{q}} $, by using (\ref{bchoi}),  we can easily get
\begin{eqnarray*}
	t\mid W_{l}^{*}\mid &\leq &\frac{r\log(\hat{\mathbf{n}})}{k_\mathbf{n}} b^{N}
	\leq  C\frac{\log(\hat{\mathbf{n}})}{k_\mathbf{n}}\hat{\mathbf{n}}^{2(1-s(1\tilde{\gamma}))/a}\\
	& \leq & C \frac{\log(\hat{\mathbf{n}})}{\hat{\mathbf{n}}^{\tilde{a} /a}},
\end{eqnarray*}
where $\tilde{a} = a\tilde{\gamma} -2(1-s(1-\tilde{\gamma}))>0$ and $\tilde{a}>0$. However, we have $t\mid W_{l}^{*}\mid<1$ for $\mathbf{n}$ large enough.
\noindent So,
$\exp\left(t W_{l}^{*}\right)\leq 1+tW_{l}^{*}+t^{2} W_{l}^{*2}$
then
$$\Ex \left(\exp\left(t W_{l}^{*}\right)\right)\leq 1+\Ex\left(t^{2} W_{l}^{*2}\right)\leq\exp\left(\Ex\left(t^{2} W_{l}^{*2}\right)\right).$$
Therefore,
$$\prod_{l=1}^{\hat{\mathbf{q}}}\Ex\left(\exp\left(tW_{l}^{*}\right)\right)\leq \exp\left(t^{2}\sum_{l=1}^{\hat{\mathbf{q}}}\Ex\left(W_{l}^{*2}\right)\right).$$
As $W_{l}^{*} $ and $W_{l}$ have the same distribution, we have
$$\sum_{l=1}^{\hat{\mathbf{q}}}\Ex\left((W_{l}^{*})^{2}\right)=Var\left(\sum_{l=1}^{\hat{\mathbf{q}}} W_{l}^{*}\right)=Var\left(\sum_{l=1}^{\hat{\mathbf{q}}}W_{l}\right)\leq S_{\mathbf{n}}+R_{\mathbf{n}}.$$
From Lemma \ref{l1}, we obtain
\begin{eqnarray*}
	\prod_{l=1}^{\hat{\mathbf{q}}}\Ex\left(\exp\left(t W_{l}^{*}\right)\right)&\leq&\exp\left(Ct^{2}k_\mathbf{n}\right)\leq
	\exp\left( Cr^{2}\frac{\log(\hat{\mathbf{n}})^{2}}{k_\mathbf{n}}\right)\longrightarrow 1, 
\end{eqnarray*}
because $\log(\hat{\mathbf{n}})^{2}/k_\mathbf{n} \rightarrow 0$ as $\mathbf{n}\to \infty.$

\noindent Then, we deduce that
\begin{eqnarray*}
	P_{12,\mathbf{n}}&\leq &C\exp\left(-t\frac{C k_\mathbf{n}(1-\sqrt{\beta})}{2^{N+1}}\right)\\
	&\leq& C\exp\left(-\frac{r C(1-\sqrt{\beta})}{2^{N+1}}\log(\hat{\mathbf{n}})\right)
	\leq C\hat{\mathbf{n}}^{-\frac{rC(1-\sqrt{\beta})}{2^{N+1}}}.
\end{eqnarray*}
Then, we have
$$N_{\varepsilon_\mathbf{n}}P_{12,\mathbf{n}}< C\hat{\mathbf{n}}^{\gamma-\tilde{\gamma}-\frac{rC(1-\sqrt{\beta})}{2^{N+1}}}.$$
Therefore, for some $r>0$ such that $\displaystyle \frac{rC(1-\sqrt{\beta})}{2^{N+1}}\tilde{\gamma}-\gamma>1$, we get
$$ \sum_{\mathbf{n}\in \mathbb{N}^{N}}N_{\varepsilon_\mathbf{n}} P_{12,\mathbf{n}}<\infty.$$
By combining the two results on $P_{11,\mathbf{n}}$ and $P_{12,\mathbf{n}}$, we get
$\displaystyle \sum_{\mathbf{n}\in \mathbb{N}^{N}} N_{\varepsilon_\mathbf{n}}P_{1,\mathbf{n}}<\infty$.\\
Using similar arguments, note that $\sum_{\mathbf{n}\in \mathbb{N}^{N}} N_{\varepsilon_\mathbf{n}}P_{2,\mathbf{n}}<\infty$.\\

\noindent Now the check of conditions $(L_2)$, $(L_3)$, $(L^{'}_{2})$ and $(L^{'}_{3})$ is based on Theorem 3.1 in \cite{dab14}. We need to show that $D_{\mathbf{n}}^{-}(\beta,x)$,  $D_{\mathbf{n}}^{+}(\beta,x)$ satisfy assumptions (H6) and (H7) used by these authors for all $(\beta,x) \in ]0,1[ \times D$. This is proved in the following lemmas where without ambiguity $D_{\mathbf{n}}$  will  denote $D_{\mathbf{n}}^{-}(\beta,x)$ or $D_{\mathbf{n}}^{+}(\beta,x)$.
\begin{lemma}\label{l2}
	Under assumption (H2) and (H6) on $\psi(.)$, we have
	$$ \hat{\mathbf{n}}D_{\mathbf{n}}^{d\theta_{0}}h_{\mathbf{n},\mathbf{s}_0}^{N\theta_{1}}\log(\hat{\mathbf{n}})^{-\theta_{2}}u_{\mathbf{n}}^{-\theta_{3}} \rightarrow \infty$$
	with
	$$\theta_{0}=\frac{s(\theta+N(d+2))}{\theta-N(s(d+4)+2d)};\qquad
	\theta_{1}=\frac{s(\theta+Nd)}{\theta-N(s(d+4)+2d)},$$
	$$\theta_{2}=\frac{s(\theta-N(d+2))}{\theta-N(s(d+4)+2d)};\qquad\theta_{3}=\frac{2(\theta+N(d+s))}{\theta-N(s(d+4)+2d)},$$
	and  $\displaystyle u_{\mathbf{n}}=\prod_{i=1}^{N}\left(\log(\log(n_{i}))\right)^{1+\varepsilon}\log(n_{i})$ for all $\varepsilon>0$.
\end{lemma}
\subsection*{Proof of Lemma \ref{l2}}

\noindent By the definition of $D_{\mathbf{n}}$ in Lemma \ref{l1}, hypotheses (H2) and (H6), we have
\begin{eqnarray*}
	\hat{\mathbf{n}}D_{\mathbf{n}}^{d\theta_{0}}h_{\mathbf{n},\mathbf{s}_{0}}^{N\theta_{1}}\log(\hat{\mathbf{n}})^{-\theta_{2}}u_{\mathbf{n}}^{-\theta_{3}}
	& \geq & C\hat{\mathbf{n}}\left(\frac{k_{\mathbf{n}}}{k_{\mathbf{n}}^{'}}\right)^{\theta_0}\left(\frac{k^{'}_\mathbf{n}}{\hat{\mathbf{n}}}\right)^{\theta_1}
	\log(\hat{\mathbf{n}})^{-\theta_2}u_{\mathbf{n}}^{-\theta_3}\\
	& \geq & C \frac{\hat{\mathbf{n}}^{1 -(\gamma -\tilde{\gamma})\theta_0-(1-\gamma)\theta_1 }}{\log(\hat{\mathbf{n}})^{\theta_2}u_{\mathbf{n}}^{\theta_3}}.
\end{eqnarray*}
Note that $ u_{\mathbf{n}}\leq \log(\tilde{n})^{N(2+\varepsilon)} \Rightarrow \frac{1}{u_{\mathbf{n}}^{\theta_3}}\geq \frac{1}{\log(\tilde{n})^{(2+\varepsilon)N\theta_3}}$, where $\displaystyle \tilde{n}=\max_{k=1,\ldots,N}n_{k}$, and
$$ \log(\hat{\mathbf{n}})\leq C\log(\tilde{n})\Rightarrow \frac{1}{\log(\hat{\mathbf{n}})^{\theta_2}}\geq C\frac{1}{\log(\tilde{n})^{\theta_2}}.$$
Since $\displaystyle \frac{n_{k}}{n_{i}}\leq C, \; \forall\; 1\leq k,i \leq N $, we deduce that $\displaystyle \hat{\mathbf{n}}\geq C \tilde{n}^{N}$ and 
\begin{eqnarray*}
	\hat{\mathbf{n}}D_{\mathbf{n}}^{d\theta_{0}}h_{\mathbf{n},\mathbf{s}_{0}}^{N\theta_{1}}\log(\hat{\mathbf{n}})^{-\theta_{2}}u_{\mathbf{n}}^{-\theta_{3}}
	& \geq &C \frac{\tilde{n}^{N\left(1 -(\gamma -\tilde{\gamma})\theta_0-(1-\gamma)\theta_1 \right)}}{\log(\tilde{n})^{\theta_2+N\theta_3(\varepsilon+2)}} \qquad  \rightarrow +\infty
\end{eqnarray*}
because 
$(1-s(1-\tilde{\gamma}))\theta>N\left ((2+s(2-\tilde{\gamma}))d+2s(2+\gamma-\tilde{\gamma})\right )$.

\begin{lemma}\label{l3}
	Under assumption (H2) and (H7) on $\psi(.)$, we have
	$$ \hat{\mathbf{n}}D_{\mathbf{n}}^{d\theta_{0}^{'}}h_{\mathbf{n},\mathbf{s}_0}^{N\theta_{1}^{'}}\log(\hat{\mathbf{n}})^{-\theta_{2}^{'}}u_{\mathbf{n}}^{-\theta_{3}^{'}} \rightarrow \infty,$$
	with
	$$\theta_{0}^{'}=\frac{s(\theta+N(d+3))}{\theta-N\left(s(d+3+2\tilde{\beta})+2(d+1)\right)};\quad \theta_{1}^{'}=\frac{s(\theta+N(d+1))}{\theta-N\left(s(d+3+2\tilde{\beta})+2(d+1)\right)}$$
	$$\theta_{2}^{'}=\frac{s(\theta-N(d+1))}{\theta-N\left(s(d+3+2\tilde{\beta})+2(d+1)\right)}; \quad \theta_{3}^{'}=\frac{2\left(\theta+N(s+d+1)\right)}{\theta-N\left(s(d+3+2\tilde{\beta})+2(d+1)\right)}.$$
\end{lemma}
\noindent The  proof of this lemma is the same as the one of Lemma~\ref{l2} and is omitted.
\subsection*{Verification of $(L_{2})$}
\noindent Let
$$f_{\mathbf{n}}\left(x,D_{\mathbf{n}}^{-}(\beta,x)\right)=\frac{1}{\hat{\mathbf{n}}h_{\mathbf{n},\mathbf{s_0}}^{N}\left(D_{\mathbf{n}}^{-}(\beta,x)\right)^{d}}\sum_{\mathbf{i}\in \mathcal{I}_\mathbf{n},\, \mathbf{i}\neq \mathbf{s_0}} K_{1}\left(\frac{x-X_{\mathbf{i}}}{D_{\mathbf{n}}^{-}(\beta,x)}\right) K_{2}\left( h_{\mathbf{n},\mathbf{s_0}}^{-1}\left\|\frac{\mathbf{s}_{0}-\mathbf{i}}{\mathbf{n}}\right\|\right),$$
and
$$f_{\mathbf{n}}\left(x,D_{\mathbf{n}}^{+}(\beta,x)\right)=\frac{1}{\hat{\mathbf{n}}h_{\mathbf{n},\mathbf{s_0}}^{N}\left(D_{\mathbf{n}}^{+}(\beta,x)\right)^{d}}\sum_{\mathbf{i}\in \mathcal{I}_\mathbf{n},\,\mathbf{i}\neq \mathbf{s_0}} K_{1}\left (\frac{x-X_{\mathbf{i}}}{D_{\mathbf{n}}^{+}(\beta,x)}\right)K_{2}\left (h_{\mathbf{n},\mathbf{s_0}}^{-1}\left \|\frac{\mathbf{s_0}-\mathbf{i}}{\mathbf{n}}\right\|\right).$$
Under the hypotheses of Lemma~\ref{th1_reg} and the results of Lemma~\ref{l2} and Lemma~\ref{l3} \cite[see][]{dab14}, we have
$$\sup_{x \in D} \left|f_{\mathbf{n}}\left(x,D_{\mathbf{n}}^{-}(\beta,x)\right)-f(x)\right|\longrightarrow 0 \qquad a.co. $$
$$\sup_{x \in D}\left|f_{\mathbf{n}}\left(x,D_{\mathbf{n}}^{+}(\beta,x)\right)-f(x)\right|\longrightarrow 0 \qquad a.co, $$
then,
$$\sup_{x \in D} \left \vert \frac{\sum_{\mathbf{i}\in \mathcal{I}_\mathbf{n},\,\mathbf{i}\neq \mathbf{s_0}}K_{1}\left(\frac{x-X_{\mathbf{i}}}{D_{\mathbf{n}}^{-}(\beta,x)}\right) K_{2}\left( h_{\mathbf{n},\mathbf{s_0}}^{-1}\left \|\frac{\mathbf{s_0}-\mathbf{i}}{\mathbf{n}}\right\|\right) }{\sum_{\mathbf{i}\in \mathcal{I}_\mathbf{n},\,\mathbf{i}\neq \mathbf{s_0}} K_{1}\left( \frac{x-X_{\mathbf{i}}}{D_{\mathbf{n}}^{+}(\beta,x)}\right) K_{2}\left (h_{\mathbf{n},\mathbf{s_0}}^{-1}\left \|\frac{\mathbf{s_0}-\mathbf{i}}{\mathbf{n}}\right\|\right)}-\beta \right \vert=\beta\sup_{x \in D} \left \vert\frac{f_{\mathbf{n}}\left(x,D_{\mathbf{n}}^{-}(\beta,x)\right)}{f_{\mathbf{n}}\left(x,D_{\mathbf{n}}^{+}(\beta,x)\right)} - 1\right \vert \rightarrow 0 \quad a.co.$$
\subsection*{Verification of $(L_{3})$}
\noindent  Under assumptions of Lemma~\ref{th1_reg}, Lemma~\ref{l2} and \ref{l3}, it follows that  (see \cite{dab14})
$$ \sup_{x \in D} \lvert c_{\mathbf{n}}\left(D_{\mathbf{n}}^{-}(\beta,x)\right)-r(x)\rvert \rightarrow 0 \quad a.co \;\mathrm{and} \; \sup_{x \in D}\lvert c_{\mathbf{n}}\left(D_{\mathbf{n}}^{+}(\beta,x)\right)-r(x)\rvert\rightarrow 0 \quad a.co.$$
\subsection*{ Proof of Lemma~\ref{th2_reg}}
\noindent The proof of this lemma is based on the results of Lemma~\ref{lemtech}. It suffices to check the conditions $(L_{2}')$ and  $(L_{3}')$. Clearly, similar arguments as those involved to prove $(L_2)$ and $(L_3)$ can be used to obtain the requested conditions.

\subsection*{Verification of $(L_{2}')$}
\noindent Under assumptions of Corollary~\ref{Cor1} and Lemmas~\ref{l2}, \ref{l3}, we have
\begin{eqnarray*}
	\sup_{x\in D}\left \vert f_{\mathbf{n}}\left(x,D_{\mathbf{n}}^{-}(\beta,x)\right)-f(x) \right\vert &=&\mathcal{O}\left(D_{\mathbf{n}}^{-}(\beta,x)\right)+\mathcal{O}\left(\left(\frac{\log(\hat{\mathbf{n}})}{\hat{\mathbf{n}}(D_{\mathbf{n}}^{-}(\beta,x))^{d}h_{\mathbf{n},\mathbf{s_0}}^{N}}\right)^{1/2}\right)\; a.co.\\
	&=&\mathcal{O}\left(\left(\frac{k_\mathbf{n}}{k^{'}_{\mathbf{n}}}\right)^{1/d}+\left(\frac{\log(\hat{\mathbf{n}})}{k_\mathbf{n}}\right)^{1/2}\right) \; a.co.,
\end{eqnarray*}
\begin{eqnarray*}
	\sup_{x\in D}\left \vert f_{\mathbf{n}}\left(x,D_{\mathbf{n}}^{+}(\beta,x)\right)-f(x)\right \vert &=&\mathcal{O}\left(D_{\mathbf{n}}^{+}(\beta,x)\right)+\mathcal{O}\left(\left(\frac{\log(\hat{\mathbf{n}})}{\hat{\mathbf{n}}(D_{\mathbf{n}}^{+}(\beta,x))^{d}h_{\mathbf{n},\mathbf{s_0}}^{N}}\right)^{1/2}\right)\;a.co.\\
	&=&\mathcal{O}\left(\left(\frac{k_\mathbf{n}}{k^{'}_{\mathbf{n}}}\right)^{1/d}+\left(\frac{\log(\hat{\mathbf{n}})}{k_\mathbf{n}}\right)^{1/2}\right). \; a.co.
\end{eqnarray*}
We conclude that
\begin{eqnarray*}
	\sup_{x\in D}\left\vert \frac{\sum_{\mathbf{i}\in \mathcal{I}_\mathbf{n},\,\mathbf{i}\neq \mathbf{s_0}}K_{1}\left(\frac{x-X_{\mathbf{i}}}{D_{\mathbf{n}}^{-}(\beta, x)}\right)K_{2}\left (h_{\mathbf{n},\mathbf{s_0}}^{-1}\left\|\frac{\mathbf{s_0}-\mathbf{i}}{\mathbf{n}}\right\|\right)}{\sum_{\mathbf{i}\in \mathcal{I}_\mathbf{n},\, \mathbf{i}\neq \mathbf{s_0}}K_{1}\left (\frac{x-X_{\mathbf{i}}}{D_{\mathbf{n}}^{+}(\beta, x)}\right) K_{2}\left (h_{\mathbf{n},\mathbf{s_0}}^{-1}\left \|\frac{\mathbf{s_0}-\mathbf{i}}{\mathbf{n}}\right\|\right)}-\beta\right \vert
	\qquad \qquad \qquad \qquad \qquad  \qquad \qquad \qquad  \\
	=\beta \sup_{x\in D}\left \vert \frac{f_{\mathbf{n}}\left(x,D_{\mathbf{n}}^{-}(\beta, x)\right)}{f_{\mathbf{n}}\left(x,D_{\mathbf{n}}^{+}(\beta, x)\right)}-1 \right \vert
	=\mathcal{O}\left(\left(\frac{k_\mathbf{n}}{k^{'}_{\mathbf{n}}}\right)^{1/d}+\left(\frac{\log(\hat{\mathbf{n}})}{k_\mathbf{n}}\right)^{1/2}\right) \, a.co.
\end{eqnarray*}
\subsection*{Verification of $(L_{3}')$}
\noindent It is relatively easy to deduce from Lemmas~\ref{l2} and \ref{l3} (\cite{dab14}) that
$$\sup_{x\in D} \left \vert c_{\mathbf{n}}\left(D_{\mathbf{n}}^{-}(\beta, x)\right)-r(x)\right \vert =\mathcal{O}\left(\left(\frac{k_\mathbf{n}}{k^{'}_{\mathbf{n}}}\right)^{1/d}+\left(\frac{\log(\hat{\mathbf{n}})}{k_\mathbf{n}}\right)^{1/2}\right) \; a.co.$$
$$
\sup_{x\in D} \left \vert c_{\mathbf{n}}\left(D_{\mathbf{n}}^{+}(\beta, x)\right)-r(x)\right \vert =\mathcal{O}\left(\left(\frac{k_\mathbf{n}}{k^{'}_{\mathbf{n}}}\right)^{1/d}+\left(\frac{\log(\hat{\mathbf{n}})}{k_\mathbf{n}}\right)^{1/2}\right) \, a.co.
$$
This yields the proof.
\subsection*{\textbf{Supplementary Materials}}
\noindent
The R code of the proposed $k$-NN predictor and classifier is available at the following \href{https://drive.google.com/file/d/1DLdXpWiiGyMfMXsQ4TVkypUHC36BYJbu/view?usp=sharing}{\textcolor{blue}{link}}. It  also allows the replicability of the environmental case study. \\

\bibliographystyle{model5-names}
\bibliography{BibkNN}
\end{document}